\documentclass[draft,12pt]{amsart}
\usepackage{amssymb,times}
\usepackage[numbers]{natbib}

\newtheorem*{Thm*}{Theorem}
\newtheorem{Thm}{Theorem}
\newtheorem{Cor}[Thm]{Corollary}
\newtheorem{Prop}[Thm]{Proposition}
\newtheorem{Lemma}{Lemma}

\theoremstyle{definition}
\newtheorem{Defn}{Definition}
\newtheorem{Notation}[Defn]{Notation}

\newtheorem{Remark}{Remark}
\newtheorem{Ex}[Remark]{Example}

\newcommand{\mf}[1]{\mathbb{#1}}
\newcommand{\mc}[1]{\mathcal{#1}}
\newcommand{\mb}[1]{\mathbf{#1}}
\newcommand{\mk}[1]{\mathfrak{#1}}

\DeclareMathOperator{\NC}{\mathit{NC}}
\DeclareMathOperator{\Falg}{\mathcal{F}_{\mathrm{alg}}}

\DeclareMathOperator{\Int}{\mathit{Int}}

\newcommand{\norm}[1]{\left\Vert#1\right\Vert}
\newcommand{\abs}[1]{\left\vert#1\right\vert}
\newcommand{\chf}[1]{\mathbf{1}_{#1}}
\newcommand{\set}[1]{\left\{#1\right\}}
\newcommand{\ip}[2]{\left \langle #1, #2 \right \rangle}
\newcommand{\state}[1]{\varphi \left[ #1 \right]}
\renewcommand{\phi}{\varphi}

\newcommand{\Cum}[1]{R \left[ #1 \right]}
\newcommand{\Span}[1]{\mathrm{Span} \left( #1 \right)}

\newcommand{\br}{\medskip\noindent}

\allowdisplaybreaks[1]

\title{Free Meixner states}
\author[M.~Anshelevich]{Michael Anshelevich}
\thanks{This work was supported in part by NSF grant DMS-0613195}
\address{Department of Mathematics, Texas A\&M University, College Station, TX 77843-3368}
\email{manshel@math.tamu.edu}
\subjclass[2000]{Primary 46L54; Secondary 05E35, 33C47}
\date{\today}

\begin{document}

\begin{abstract}
Free Meixner states are a class of functionals on non-commutative polynomials introduced in \cite{AnsMulti-Sheffer}. They are characterized by a resolvent-type form for the generating function of their orthogonal polynomials, by a recursion relation for those polynomials, or by a second-order non-commutative differential equation satisfied by their free cumulant functional. In this paper, we construct an operator model for free Meixner states. By combinatorial methods, we also derive an operator model for their free cumulant functionals. This, in turn, allows us to construct a number of examples. Many of these examples are shown to be trivial, in the sense of being free products of functionals which depend on only a single variable, or rotations of such free products. On the other hand, the multinomial distribution is a free Meixner state and is not a product. Neither is a large class of tracial free Meixner states which are analogous to the simple quadratic exponential families in statistics.
\end{abstract}

\maketitle

\section{Introduction}
\noindent
The subject of this paper are states and orthogonal polynomials in non-commuting variables. The definition is straightforward. The usual orthogonal polynomials are obtained by starting with a measure $\mu$ on $\mf{R}^d$, thinking of $\mf{R}[x_1, x_2, \ldots, x_d]$ as a vector space with the (pre-)inner product
\[
\ip{P}{Q} = \int_{\mf{R}^d} P(\mb{x}) Q(\mb{x}) \,d\mu(\mb{x}),
\]
and applying the Gram-Schmidt procedure to the monomials $\set{x_{u(1)} x_{u(2)} \ldots x_{u(n)}}$. In the non-commutative case, one starts directly with a positive linear functional (state) $\phi$ on the algebra of non-commutative polynomials $\mf{R} \langle x_1, x_2, \ldots, x_n \rangle$, and orthogonalizes the monomials in non-commuting variables with respect to the inner product
\[
\ip{P}{Q} = \state{P^\ast(\mb{x}) Q(\mb{x})}.
\]

\br
Among the general ``non-commutative measures'' and polynomials orthogonal with respect to them, there is a specific class of what is appropriate to call \emph{free Meixner states}. The classical Meixner class \cite{Meixner} consists of familiar distributions---normal, Poisson, gamma, negative binomial, Meix\-ner, and binomial---which, somewhat less familiarly, share a number of common properties: their orthogonal polynomials have exponential-form generating functions, they satisfy a quadratic regression property \cite{Laha-Lukacs}, they generate quadratic natural exponential families \cite{Morris}, they are quadratic harnesses \cite{Wes-commutative}, they are induced by representations of $\mk{su}(1,1)$ \cite{Koelink-Convolutions}, and they have explicit linearization coefficient formulas \cite{KimZeng}. The multivariate Meixner distributions have also been investigated, frequently in the guise of quadratic exponential families \cite{Casalis-Simple-quadratic,Pommeret-Test}, though a complete classification is still lacking. Even the infinite-dimensional case was considered \cite{Sniady-SWN,Lytvynov-Meixner}.

\br
In \cite{AnsMeixner}, I introduced the free Meixner polynomials, which are a family of orthogonal polynomials in one variable. The term ``free'' refers to their relation to free probability, see \cite{VDN,Nica-Speicher-book} for an introduction. As a matter of fact, these polynomials have been found independently both before and after my work, for example in \cite{Sze22,CTConstant,Freeman,SaiConstant,Kubo-IDAQP}. They share a number of the Meixner properties listed above, as long they are properly translated into the ``free'' context, see my original paper and also \cite{Boz-Bryc}. Some of the corresponding distributions also appear in random matrix theory, as the limiting distributions in the Gaussian, Wishart, and Jacobi ensembles.

\br
In \cite{AnsMulti-Sheffer} I started the investigation of multivariate free Meixner distributions, which are states on the algebra of non-commutative polynomials. I continue their study in Section~\ref{Section:Meixner}. The main new tool is to represent these states as joint distributions of certain operators on a Fock space, following the more general construction in~\cite{AnsMonic}. I use this machinery, in combination with combinatorial methods, to find explicit formulas for the free cumulants of these states. This provides an explanation for the one-variable results in Section 3.1 of \cite{AnsMeixner} and Proposition 2.2 of \cite{Boz-Bryc}, and is the first main result of the paper. The operator representation of the state also allows me to handle states that are not necessarily faithful, thus answering a question of the referee of \cite{AnsMulti-Sheffer}, where only faithful free Meixner states were considered.

\br
Having an explicit representation for the cumulants, and being able to handle non-faithful states, allows me to describe a number of examples, which is done in Section~\ref{Section:Examples}. Among the usual multivariate Meixner distributions, two are familiar, namely the multivariate normal and the multinomial distributions. It is well known that the free analog of the multivariate normal distribution is the distribution of a free semicircular system, see Section~\ref{Subsec:Semicircular}. The second question treated in this paper is: what is the ``free'' multinomial distribution? I show that the basic multinomial distribution \emph{itself} also belongs to the free Meixner class. In particular, this allows me to calculate the distribution of a free sum of $d$-tuples of orthogonal projections.

\br
Among states on non-commutative algebras, \emph{traces} form an important class. The final result in this paper provides a way to construct a large family of non-trivial, tracial free Meixner states. These turn out to be analogs of simple quadratic exponential families.

\section{Preliminaries}
\noindent
Variables in this paper will typically come in $d$-tuples, which will be denoted using the bold font: $\mb{x} = (x_1, x_2, \ldots, x_d)$, and the same for $\mb{z}, \mb{S}$, etc.

\subsection{Polynomials}
Let $\mf{R}\langle \mb{x} \rangle = \mf{R}\langle x_1, x_2, \ldots, x_d \rangle$ be all the polynomials with real coefficients in $d$ non-commuting variables. \emph{Multi-indices} are elements $\vec{u} \in \set{1, \ldots, d}^k$ for $k \geq 0$; for $\abs{\vec{u}} = 0$ denote $\vec{u}$ by $\emptyset$. Monomials in non-commuting variables $(x_1, \ldots, x_d)$ are indexed by such multi-indices:
\[
x_{\vec{u}} = x_{u(1)} \ldots x_{u(k)}.
\]
Note that our use of the term ``multi-index'' is different from the usual one, which is more suited for indexing monomials in commuting variables.

\br
For two multi-indices $\vec{u}, \vec{v}$, denote by $(\vec{u}, \vec{v})$ their concatenation. For $\vec{u}$ with $\abs{\vec{u}} = k$, denote
\[
(\vec{u})^{op} = (u(k), \ldots, u(2), u(1)).
\]
Define an involution on $\mf{R}\langle \mb{x} \rangle$ via the $\mf{R}$-linear extension of
\[
(x_{\vec{u}})^\ast = x_{(\vec{u})^{op}}.
\]

\br
A \emph{monic polynomial family} in $\mb{x}$ is a family $\set{P_{\vec{u}}(\mb{x})}$ indexed by all multi-indices
\[
\bigcup_{k=1}^\infty \set{\vec{u} \in \set{1, \ldots, d}^k}
\]
(with $P_{\emptyset} = 1$ being understood) such that\[
P_{\vec{u}}(\mb{x}) = x_{\vec{u}} + \textsl{lower-order terms}.
\]
Note that $P_{\vec{u}}^\ast \neq P_{(\vec{u})^{op}}$ in general.

\begin{Defn}
\label{Defn:State}
A \emph{state} on $\mf{R} \langle \mb{x} \rangle$ is a functional
\[
\phi: \mf{R} \langle x_1, x_2, \ldots, x_d \rangle \rightarrow \mf{R}
\]
that is linear, compatible with the $\ast$-operation, that is for any $P$,
\[
\state{P} = \state{P^\ast},
\]
unital, that is $\state{1} = 1$, and positive, that is for any $P$,
\[
\state{P^\ast P} \geq 0.
\]
A state is \emph{faithful} if in the preceding equation, the equality holds only for $P = 0$. Unless noted otherwise, the states in this paper are \emph{not} assumed to be faithful.

\br
The numbers $\state{x_{\vec{u}}}$ are called the \emph{moments} of $\phi$.

\br
A state $\phi$ induces the pre-inner product
\[
\ip{P}{Q}_\phi = \state{P^\ast Q} = \ip{Q}{P}_\phi
\]
and the seminorm
\[
\norm{P}_\phi = \sqrt{\state{P^\ast P}}.
\]
Throughout the paper, we will typically drop $\phi$ from the notation, and denote the inner product and norm it induces simply by $\ip{\cdot}{\cdot}$, $\norm{\cdot}$.

\br
We may think of $\phi$ is a ``joint distribution'' of ``random variables'' $(x_1, x_2, \ldots, x_d)$. In the remainder of the paper, as we did in \cite{AnsMulti-Sheffer}, we will assume that under the state $\phi$, the variables have zero mean and identity covariance,
\[
\state{x_i} = 0, \qquad \state{x_i x_j} = \delta_{ij}.
\]
The last assumption is made primarily so that equation~\eqref{PDE} has a clean form. In Section~\ref{Subsec:Covariance} we briefly describe how to modify the results if that assumption is dropped.
\end{Defn}

\subsection{Monic orthogonal polynomials states}
\label{Subsec:MOPS}
\begin{Defn}
A state has a \emph{monic orthogonal polynomial system}, or MOPS, if for any multi-index $\vec{u}$, there is a monic polynomial $P_{\vec{u}}$ with leading term $x_{\vec{u}}$, such that these polynomials are orthogonal with respect to $\phi$, that is,
\[
\ip{P_{\vec{u}}}{P_{\vec{v}}} = 0
\]
for $\vec{u} \neq \vec{v}$.
\end{Defn}

\noindent
Note that the same abbreviation is used in~\cite{Dumitriu-MOPS} to denote a class of multivariate orthogonal polynomials systems, which is different from ours.

\br
States that have MOPS were characterized in \cite{AnsMonic}. We briefly summarize the results of that paper which we will use in the next section.

\subsubsection{Fock space construction I}
\label{Subsubsec:General-Fock}
Let $\mc{H} = \mf{C}^d$, with the canonical orthonormal basis $e_1, e_2, \ldots, e_d$. Define the (algebraic) full Fock space of $\mc{H}$ to be
\[
\Falg(\mc{H}) = \bigoplus_{k=0}^\infty \mc{H}^{\otimes k}
\]
Equivalently, $\Falg(\mc{H})$ is the vector space of non-commutative polynomials in $e_1, e_2, \ldots, e_d$. Following convention, we will denote the generating vector in $\mc{H}^{\otimes 0} = \mf{C}$ by $\Omega$ instead of $1$.

\br
For $i = 1, 2, \ldots, d$, define $a_i^+$ and $a_i^-$ to be the usual (left) free creation and annihilation operators,
\begin{align*}
a_i^+ & \left(e_{u(1)} \otimes e_{u(2)} \otimes \ldots \otimes e_{u(k)} \right) = e_i \otimes e_{u(1)} \otimes e_{u(2)} \otimes \ldots \otimes e_{u(k)}, \\
a_i^-      & (e_j) = \ip{e_i}{e_j} \Omega = \delta_{i j} \Omega, \\
a_i^-      & \left(e_{u(1)} \otimes e_{u(2)} \otimes \ldots \otimes e_{u(k)} \right) = \ip{e_i}{e_{u(1)}} e_{u(2)} \otimes \ldots \otimes e_{u(k)}.
\end{align*}

\br
For each $k \geq 2$ let $\mc{C}^{(k)}$ be an operator
\[
\mc{C}^{(k)}: \mc{H}^{\otimes k} \rightarrow \mc{H}^{\otimes k}.
\]
We think of each $\mc{C}^{(k)}$ as a $d^k \times d^k$ matrix. Assume that for each $k$, $\mc{C}^{(k)}$ is diagonal and $\mc{C}^{(k)} \geq 0$. It is convenient to also take $\mc{C}^{(1)} = I$; this corresponds to the identity covariance. Similarly, for each $i = 1, 2, \ldots, d$ and each $k \geq 1$, let $\mc{T}_i^{(k)}$ be an operator
\[
\mc{T}_i^{(k)}: \mc{H}^{\otimes k} \rightarrow \mc{H}^{\otimes k}.
\]
Assume that $\mc{T}_i^{(k)}$ and $\mc{C}^{(j)}$ satisfy a commutation relation (see \cite{AnsMonic}). We will denote by $\mc{T}_i$ and $\mc{C}$ the operators acting as $\mc{T}_i^{(k)}$ and $\mc{C}^{(k)}$ on each component. Finally, let $\tilde{a}_i^- = a_i^- \mc{C}$ and
\[
\mc{X}_i = a_i^+ + \mc{T}_i + \tilde{a}_i^-.
\]
With the appropriate choice of the inner product $\ip{\cdot}{\cdot}_{\mc{C}}$ on the completion $\mc{F}_{\mc{C}}(\mc{H})$ of the quotient of $\Falg(\mc{H})$, all the operators $a_i^+, \mc{T}_i, \tilde{a}_i^-$ factor through to $\mc{F}_{\mc{C}}(\mc{H})$, and each $\mc{X}_i$ is a symmetric operator on it.

\begin{Thm}(Part of Theorem 2 of \cite{AnsMonic})
\label{Thm:Monic-states}
Let $\phi$ be a state on $\mf{R} \langle \mb{x} \rangle$. The following are equivalent:
\begin{enumerate}
\item
The state $\phi$ has a monic orthogonal polynomial system.
\item
There is a family of polynomials $\set{P_{\vec{u}}}$ such that $\state{P_{\vec{u}}} = 0$ for all $\vec{u} \neq \emptyset$ and they satisfy a recursion relation
\begin{align*}
x_i & = P_i + B_{i, \emptyset, \emptyset}, \\
x_i P_u & = P_{(i, u)} + \sum_{w=1}^d B_{i, w, u} P_{w} + \delta_{i, u} C_u, \\
x_i P_{\vec{u}} & = P_{(i, \vec{u})} + \sum_{\abs{\vec{w}} = \abs{\vec{u}}} B_{i, \vec{w}, \vec{u}} P_{\vec{w}} + \delta_{i, u(1)} C_{\vec{u}} P_{(u(2), u(3), \ldots, u(k))},
\end{align*}
with $C_{\vec{u}} \geq 0$ and, denoting $\vec{s}_j = (s(j), \ldots, s(k))$,
\[
B_{i, \vec{s}, \vec{u}} \prod_{j=1}^k C_{\vec{s}_j} = B_{i, \vec{u}, \vec{s}} \prod_{j=1}^k C_{\vec{u}_j}.
\]
\item

For some choice of the matrices $\mc{C}^{(k)}$ and $\mc{T}_i^{(k)}$ as in Section~\ref{Subsubsec:General-Fock}, the state $\phi$ has a Fock space representation $\phi_{\mc{C}, \set{\mc{T}_i}}$ as
\begin{equation*}
\state{P(x_1, x_2, \ldots, x_d)} = \ip{\Omega}{P(\mc{X}_1, \mc{X}_2, \ldots, \mc{X}_d) \Omega}.
\end{equation*}
\end{enumerate}
\end{Thm}

\br
We will also need the following relation between the operators in part (c) and coefficients in part (b) of the theorem:
\begin{equation}
\label{Expansion-T}
\mc{T}_i(e_{u(1)} \otimes \ldots \otimes e_{u(k)}) = \sum_{\abs{\vec{w}} = k} B_{i, \vec{w}, \vec{u}} e_{w(1)} \otimes \ldots \otimes e_{w(k)}
\end{equation}
and
\begin{equation}
\label{Expansion-C}
\mc{C}(e_{u(1)} \otimes \ldots \otimes e_{u(k)}) = C_{\vec{u}} e_{u(1)} \otimes \ldots \otimes e_{u(k)}.
\end{equation}

\subsection{Fock space construction II}
\label{Subsec:Fock2}
The following construction is a particular case of the construction in Section~\ref{Subsubsec:General-Fock}, but this time we provide full details. As before, let $\mc{H} = \mf{C}^d$, with the canonical basis $e_1, e_2, \ldots, e_d$, denote its  (algebraic) full Fock space by $\Falg(\mc{H})$, and the generator of the zeroth component by $\Omega$. Let $C$ be an n operator on $\mc{H} \otimes \mc{H}$, which we identify with its $d^2 \times d^2$ matrix in the standard basis. Assume that $C$ is diagonal and
\begin{equation}
\label{C-positive}
(I \otimes I) + C \geq 0,
\end{equation}
where $I$ will always denote the identity operator on $\mc{H}$. On $\Falg(\mc{H})$, define a new inner product using the non-negative kernel
\[
K_C = \bigl(I^{\otimes (k-2)} \otimes (I^{\otimes 2} + C)\bigr) \ldots \bigl(I \otimes (I^{\otimes 2} + C) \otimes I^{\otimes (k-3)}\bigr) \bigl((I^{\otimes 2} + C) \otimes I^{\otimes (k-2)}\bigr)
\]
on each $\mc{H}^{\otimes k}$, and denote the completion of $\Falg(\mc{H})$ with respect to this inner product $\mc{F}_C(\mc{H})$. If the inner product is degenerate, first factor out the subspace of vectors of length zero, and then complete.

\br
For $i = 1, 2, \ldots, d$, let $a_i^+$ and $a_i^-$ be the usual (left) free creation and annihilation operators as defined in Section~\ref{Subsubsec:General-Fock}. Let $T_1, \ldots, T_d$ be operators on $\mc{H}$ which we identify with their $d \times d$ matrices. Assume that each $T_i$ is symmetric and
\[
(T_i \otimes I) C = C (T_i \otimes I).
\]
With a slight abuse of notation, we will denote
\[
T_i = T_i \otimes I^{\otimes (k-1)} \text{ on } \mc{H}^{\otimes k}
\]
and
\[
\tilde{a}_i = a_i^- (C \otimes I^{\otimes (k-2)}) \text{ on } \mc{H}^{\otimes k}.
\]
Note that
\begin{equation}
\label{Zero}
a_i^- \Omega = T_i \Omega = \tilde{a}_i \Omega = 0 \text{ and } \tilde{a}_i = 0 \text{ on } \mc{H}.
\end{equation}
It follows from the general construction in Section~\ref{Subsubsec:General-Fock} that all the operators
\[
X_i = a_i^+ + a_i^- + T_i + \tilde{a}_i
\]
factor through to $\mc{F}_C(\mc{H})$.

\begin{Defn}
\label{Defn:Fock-state}
The Fock state $\phi = \phi_{C, \set{T_i}}$ on $\mf{R} \langle \mb{x} \rangle$ determined by such $C$ and $T_i$  is the state
\begin{equation*}
\state{P(x_1, x_2, \ldots, x_d)} = \ip{\Omega}{P(X_1, X_2, \ldots, X_d) \Omega} = \ip{\Omega}{P(X_1, X_2, \ldots, X_d) \Omega}_C.
\end{equation*}
\end{Defn}

\subsection{Non-crossing partitions}
A \emph{partition} $\pi$ of a set $V \subset \mf{Z}$ is a collection of disjoint subsets of $V$ (classes of $\pi$), $\pi = (B_1, B_2, \ldots, B_k)$, whose union equals $V$. Most of the time we will be interested in partitions of $\set{1, 2, \ldots, n}$. Partitions form a partially ordered set (in fact a lattice) under the operation of refinement, so that the largest partition is $\hat{1} = \bigl( \set{1, 2, \ldots, n} \bigr)$ and the smallest partition is $\hat{0} = \bigl( \set{1}, \set{2}, \ldots, \set{n} \bigr)$. We will use $i \stackrel{\pi}{\sim} j$ to denote that $i, j$ lie in the same class of $\pi$.

\br
Let $\NC(V)$ denote the collection of non-crossing partitions of $V$, which are partitions $\pi$ such that
\[
i \stackrel{\pi}{\sim} i', j \stackrel{\pi}{\sim} j', i \stackrel{\pi}{\not \sim} j, i < j < i' \Rightarrow i < j' < i'.
\]
Equivalently, a partition is non-crossing if and only if one of its classes is an interval and the restriction of the partition to the complement of this class is non-crossing. Non-crossing partitions are a sub-lattice of the lattice of all partitions. For each $n$, let $\NC(n)$ denote the lattice of non-crossing partitions of the set $\set{1, 2, \ldots, n}$. We will also denote by $\NC_0(V)$ all non-crossing partitions with no singletons (one-element classes), and by $\NC'(V)$ all the partitions $\pi$ such that
\[
\min V \stackrel{\pi}{\sim} \max V.
\]
Equivalently, partitions in $\NC'(V)$ have a single outer class---the one that contains both $\min V$ and $\max V$---in the terminology of \cite{BLS96}. (A class $B \in \pi$ is outer if there do \emph{not} exist $i, i' \not \in B$, $j \in B$ with $i \stackrel{\pi}{\sim} i'$ and $i < j < i'$.) See \cite{Nica-Speicher-book} or \cite{Stanley-volume-1} for more details on the relevant combinatorics.

\subsection{Free cumulants}
The free cumulant functional $R$ corresponding to a state $\phi$ is the linear functional on $\mf{R} \langle \mb{x} \rangle$ defined recursively by $\Cum{1} = 0$ and for $\abs{\vec{u}} = n$,
\begin{equation}
\label{Cumulants-definition}
\Cum{x_{\vec{u}}}
= \state{x_{\vec{u}}} - \sum_{\substack{\pi \in \NC(n), \\ \pi \neq \hat{1}}} \prod_{B \in \pi} \Cum{\prod_{i \in B} x_{u(i)}},
\end{equation}
which expresses $\Cum{x_{\vec{u}}}$ in terms of the joint moments and sums of products of lower-order free cumulants. From these, we can form the free cumulant generating function of $\phi$ via
\begin{equation}
\label{Non-crossing}
R(z_1, z_2, \ldots, z_d)
= \sum_{n=1}^\infty \sum_{\abs{\vec{u}} = n} \Cum{x_{\vec{u}}} z_{\vec{u}},
\end{equation}
where $\mb{z} = (z_1, \ldots, z_d)$ are non-commuting indeterminates. One can also define $R$ using an implicit functional relation involving the moment generating function of $\phi$, see Corollary~16.16 of \cite{Nica-Speicher-book}.

\subsection{Words and partitions}
\label{Subsec:Facts}
In this section, we collect a number of facts that will be useful in the proof of the next two theorems. Note that in many places, operators are considered as acting on $\Falg(\mc{H})$, with a degenerate inner product, rather than on $\mc{F}_C(\mc{H})$.

\begin{Lemma}
Let $\vec{u}$ be a multi-index indexed by a set $V \subset \mf{Z}$, and $W = \prod_{i \in V} W(i)$ be a word with $W(i)$ equal to $a^+_{u(i)}, T_{u(i)}, a_{u(i)}^-$, or $\tilde{a}_{u(i)}$.  If
\[
\ip{\Omega}{\prod_{i \in V} W(i) \Omega} \neq 0,
\]
then
\begin{equation}
\label{Catalan-walk}
\begin{split}
W(\min V) = a_{u(\min V)}^-, \quad W(\max V) = a_{u(\max V)}^+, & \\
\forall i \in V, \abs{\set{j \in V | j \geq i, W(j) = a_{u(j)}^- \text{ or } W(j) = \tilde{a}_{u(j)}}} & \leq \abs{\set{j \in V | j \geq i, W(j) = a_{u(j)}^+}}, \\
\abs{\set{j \in V | W(j) = a_{u(j)}^- \text{ or } W(j) = \tilde{a}_{u(j)}}} & = \abs{\set{j \in V | W(j) = a_{u(j)}^+}},
\end{split}
\end{equation}
and
\begin{multline}
\label{Level-two}
\abs{\set{j \in V | j \geq i, W(j) = a_{u(j)}^- \text{ or } W(j) = \tilde{a}_{u(j)}}} = \abs{\set{j \in V | j \geq i, W(j) = a_{u(j)}^+}} \\
\Rightarrow W(i) = a_{u(i)}^-.
\end{multline}
\end{Lemma}

\begin{proof}
This follows from the fact that if $\eta \in \mc{H}^{\otimes k}$, then $a_i^+(\eta) \in \mc{H}^{\otimes (k+1)}$, $T_i(\eta) \in \mc{H}^{\otimes k}$, and $a_i^-(\eta), \tilde{a}_i(\eta) \in \mc{H}^{\otimes (k-1)}$, and equation~\eqref{Zero}.
\end{proof}

\noindent
In combinatorics, equation~\eqref{Catalan-walk} is related to the notion of a Motzkin path. More generally, our operator representations are closely related to a common way of representing moments as sums over lattice paths \cite{Flajolet,Viennot-Short}, but in the multivariate case we find the operator formulation more useful.

\begin{Notation}
Let
\[
\mc{W}_n(\vec{u}) = \bigl\{W = W(1) W(2) \ldots W(n) \text{ satisfying conditions } \eqref{Catalan-walk} \text{ and } \eqref{Level-two} \text{ for } V = \set{1, \ldots, n}\bigr\},
\]
and for a general subset $V \subset \mf{Z}$, define $\mc{W}_V(\vec{u})$ similarly.
\end{Notation}

\begin{Lemma}
\label{Lemma:Bijection}
For any multi-index $\vec{u}$, partition $\pi \in \NC_0(n)$, $\pi = (V_1, V_2, \ldots, V_k)$ and partitions $\sigma_j \in \NC_0'(V_j)$, $j = 1, 2, \ldots, k$, define a word $W = \beta_{\vec{u}}(\pi; \sigma_1, \ldots, \sigma_k)$ by
\begin{equation}
\label{Partition-word}
W(i) =
\begin{cases}
a_{u(i)}^+, & i \in B \in \sigma_j, i = \max B, \\
a_{u(i)}^-, & i \in V_j, i = \min V_j, \\
\tilde{a}_{u(i)}, & i \in B \in \sigma_j, i = \min B, i \neq \min V_j, \\
T_{u(i)}, & \text{otherwise}.
\end{cases}
\end{equation}
Then $W \in \mc{W}_n(\vec{u})$, and for each $V \in \pi$, $W$ restricted to $V$ is in $\mc{W}_V(\vec{u}:V)$, where $(\vec{u}:V)$ is the sub-multi-index of $\vec{u}$ indexed by the elements of $V$. Moreover, for each $\vec{u}$, $\beta_{\vec{u}}$ is a bijection.
\end{Lemma}

\begin{proof}
Let $W = \beta_{\vec{u}}(\pi; \sigma_1, \ldots, \sigma_k)$. Condition~\eqref{Catalan-walk} for the whole set $\set{1, 2, \ldots, n}$ (respectively, for $V_j$) follows from the definition of $\beta$ and the fact that $\pi, \sigma_1, \ldots, \sigma_k$ (respectively, $\sigma_j$) are non-crossing. Condition~\eqref{Level-two} follows from the definition that the minima of the outer classes of the partition $\pi$ (respectively, $\sigma_j$) are all $a^-$.

\br
Conversely, let $W \in \mc{W}_n(\vec{u})$. Let $\Lambda \subset \set{1, 2, \ldots, n}$,
\[
\Lambda = \set{j | W(j) \neq T_{u(j)}}.
\]
It follows from Proposition~2.13 and Exercise~8.23 of \cite{Nica-Speicher-book} that, as long as $W$ restricted to $\Lambda$ satisfies condition~\eqref{Catalan-walk}, there is a unique non-crossing pair partition $\pi' \in \NC(\Lambda)$ such that for any $B \in \pi'$,
\begin{align*}
i = \min B & \Leftrightarrow W(i) = a_{u(i)}^- \text{ or } \tilde{a}_{u(i)}, \\
i' = \max B & \Leftrightarrow W(i') = a_{u(i')}^+.
\end{align*}
Note that $W(1) = a_{u(1)}^-$, so $(1, i') \in \pi'$ for some $i' > 1$. Moreover, $W(i'+1) \ldots W(n) \Omega \in \mc{H}^{\otimes 0}$, so by condition~\eqref{Level-two}, $W(i'+1) = a_{u(i'+1)}^-$. Thus $(i'+1, j') \in \pi'$ for some $j'$, etc. ending with $(s,n) \in \pi'$. It follows that for any $j \in \set{1, \ldots, n}$, there exist $i \stackrel{\pi'}{\sim} i'$ such that $i \leq j \leq i'$ and $W(i) = a_{u(i)}^-$. For each $j$, choose the largest $i$ such that $i \stackrel{\pi'}{\sim} i'$, $i \leq j \leq i'$, and $W(i) = a_{u(i)}^-$, and require that $i \stackrel{\pi}{\sim} j \stackrel{\pi}{\sim} i'$. Similarly, for each class $V_s \in \pi$ and each $j \in V_s$, choose the largest $i \in V_s$ such that $i \stackrel{\pi'}{\sim} i'$ and $i \leq j \leq i'$, and require that $i \stackrel{\sigma_s}{\sim} j \stackrel{\sigma_s}{\sim} i'$. Pictorially, we draw the integers $1, 2, \ldots, n$ on a line, draw the pair classes or $\pi'$ as arcs connecting each $i$ with the corresponding $i'$ above the line, and then connect each of the other elements to the arc immediately above it.
\end{proof}

\begin{Lemma}
\label{Lemma:Factor}
For any $W \in \mc{W}_n(\vec{u})$, let $(\pi; \sigma_1, \ldots, \sigma_k) = \beta_{\vec{u}}^{-1}(W)$ with $\pi = (V_1, V_2, \ldots, V_k)$. Then
\[
\ip{\Omega}{W(1) W(2) \ldots W(n) \Omega} = \prod_{i=1}^k \ip{\Omega}{\prod_{j \in V_i} W(j) \Omega}.
\]
\end{Lemma}

\begin{proof}
Since $\pi$ is a non-crossing partition, it has a class $V$ that is an interval,
\[
V = [i, i'] = \set{j | i \leq j \leq i'}.
\]
Since $\pi$ restricted to $\set{1, \ldots, n} \backslash V$ is still a non-crossing partition, it suffices to show that
\[
\ip{\Omega}{W(1) W(2) \ldots W(n) \Omega} = \ip{\Omega}{\prod_{j=1}^{i-1} W(j) \prod_{j=i'+1}^n W(j) \Omega} \ip{\Omega}{\prod_{j=i}^{i'} W(j) \Omega}.
\]
Denote
\[
\eta = W(i'+1) \ldots W(n) \Omega \in \mc{H}^{\otimes m}.
\]
We now show that for any $i < j \leq i'$,
\[
W(j) \ldots W(n) \Omega = \zeta_j \otimes \eta
\]
for
\[
\zeta_j = W(j) \ldots W(i') \Omega.
\]
The proof is by induction.
\[
W(i') \eta = a_{u(i')}^+ \eta = e_{u(i')} \otimes \eta = (W(i') \Omega) \otimes \eta.
\]
If $W(j) = a_{u(j)}^+$, $\zeta_j = e_{u(j)} \otimes \zeta_{j+1}$. If $W(j) = T_{u(j)}$, then $\zeta_j = T_{u(j)} \zeta_{j+1}$. $W(j)$ cannot equal $a_{u(j)}^-$. Finally, it follows from condition~\eqref{Level-two} applied to $V = [i, i']$ that for all $j$, $i < j \leq i'$,
\[
W(j+1) \ldots W(n) \Omega \in \mc{H}^{\otimes s}
\]
with $s > m$. Thus $W(j)$ may equal $\tilde{a}_{u(j)}$ only if $s \geq m + 2$, otherwise
\[
W(j) W(j+1) \ldots W(n) \Omega \in \mc{H}^{\otimes m}.
\]
But if $s \geq m + 2$, $\zeta_j = a_{u(j)}^- C \zeta_{j+1}$.

\br
It follows that also
\[
W(i) \ldots W(n) \Omega = (W(i) \ldots W(i') \Omega) \otimes \eta = \ip{\Omega}{W(i) \ldots W(i') \Omega} \eta.
\]
Thus
\[
\begin{split}
\ip{\Omega}{W(1) W(2) \ldots W(n) \Omega}
& = \ip{\Omega}{W(1) \ldots W(i-1) \ip{\Omega}{W(i) \ldots W(i') \Omega} \eta} \\
& = \ip{\Omega}{W(1) \ldots W(i-1) \eta} \ip{\Omega}{W(i) \ldots W(i') \Omega} \\
& = \ip{\Omega}{\prod_{j=1}^{i-1} W(j) \prod_{j=i'+1}^n W(j) \Omega} \ip{\Omega}{\prod_{j=i}^{i'} W(j) \Omega}.
\qedhere
\end{split}
\]
\end{proof}

\begin{Notation}
\label{Notation:Covered-bijection}
For $V \subset \mf{Z}$, denote
\[
\begin{split}
\mc{W}_V'(\vec{u}) = \{W \in \mc{W}_V(\vec{u}) |& W(\min V) = a_{u(\min V)}^-, W(\max V) = a_{u(\max V)}^+, \\
&\quad \text{ and none of the other $W(i)$ are equal to } a_{u(i)}^-\}.
\end{split}
\]
The partition $\pi$ corresponding to any such $W$ has only one class, $\pi = (V) \in \NC(V)$, and
\[
\beta_{\vec{u}}^{-1}(\mc{W}_V'(\vec{u})) = \set{((V), \sigma) | \sigma \in \NC_0'(V)} \cong \NC_0'(V).
\]
Denote
\[
\Theta(\sigma; V, \vec{u}) = \ip{\Omega}{\beta_{\vec{u}}((V), \sigma) \Omega}.
\]
\end{Notation}

\begin{Lemma}
If $W \in \mc{W}_n(\vec{u})$ and $\beta_{\vec{u}}^{-1}(W) = (\pi, \sigma_1, \ldots, \sigma_k)$, $\pi = (V_1, V_2, \ldots, V_k)$, then
\[
\ip{\Omega}{W(1) \ldots W(n) \Omega}
= \prod_{j=1}^k \ip{\Omega}{\prod_{i \in V_j} W(i) \Omega}
= \prod_{j=1}^k \Theta(\sigma_j; V_j, (\vec{u}:V_j)),
\]
where $(\vec{u}:V_j)$ is the sub-multi-index of $\vec{u}$ indexed by the elements of $V_j$.
\end{Lemma}

\begin{proof}
This follows from Lemma~\ref{Lemma:Factor} using Notation~\ref{Notation:Covered-bijection}.
\end{proof}

\begin{Lemma}
\label{Lemma:Factor2}
Suppose that $W \in \mc{W}_n'(\vec{u})$ such that $W(1) = a_{u(1)}^- = a_{j}^-$ and $W(2) = \tilde{a}_{u(2)} = \tilde{a}_{i}$. Then $\beta_{\vec{u}}^{-1}(W) = ((\set{1, \ldots, n}), \sigma)$. It follows from condition~\eqref{Partition-word} that for $2 \in B \in \sigma$, we have $2 = \min B$. Let $k = \max B$. Then $W(k) = a_{u(k)}^+$ and
\begin{multline*}
\ip{\Omega}{W(1) W(2) \ldots W(k) \ldots W(n) \Omega}
= \ip{\Omega}{a_{j}^- \tilde{a}_{i} W(3) \ldots a_{u(k)}^+ W(k+1) \ldots W(n-1) a_{u(n)}^+ \Omega} \\
= C_{ij} \ip{\Omega}{a_j^- W(k+1) \ldots W(n-1) a_{u(n)}^+ \Omega} \ip{\Omega}{a_i^- W(3) \ldots a_{u(k)}^+ \Omega}.
\end{multline*}
Moreover, the map
\[
\begin{split}
\{W \in \mc{W}_n'(\vec{u}) | & W(1) = a_{u(1)}^- = a_{j}^-, W(2) = \tilde{a}_{u(2)} = \tilde{a}_{i}\} \\
& \rightarrow \bigcup_{k=3}^{n-1} \mc{W}_{\set{1, k+1, \ldots, n}}'\bigl((\vec{u}:\set{1, k+1, \ldots, n})\bigr) \times \mc{W}_{\set{2, \ldots, k}}'\bigl((\vec{u}:\set{2, \ldots, k})\bigr) \\
& \cong \bigcup_{k=3}^{n-1} \NC_0'(\set{1, k+1, \ldots, n}) \times \NC_0'(\set{2, \ldots, k})
\end{split}
\]
is a bijection.
\end{Lemma}

\begin{proof}
By the same method as in Lemma~\ref{Lemma:Factor}, we deduce that
\[
W(3) \ldots a_{u(k)}^+ W(k+1) \ldots W(n-1) a_{u(n)}^+ \Omega
= \bigl(W(3) \ldots a_{u(k)}^+ \Omega \bigr) \otimes \bigl(W(k+1) \ldots W(n-1) a_{u(n)}^+ \Omega \bigr).
\]
The inner product of this vector with $C(e_i \otimes e_j)$ is the desired expression.
\end{proof}

\begin{Lemma}
\label{Lemma:Last}
Suppose that $C_{ij} = C(e_i \otimes e_j) = c$ for all $i, j$. Let $\sigma \in NC_0'(n)$,
\[
\sigma = \Bigl( \set{b_{1,1}, \ldots, b_{1, j(1)}}, \ldots, \set{b_{1, k}, \ldots, b_{k, j(k)}} \Bigr),
\]
where each class is ordered and $b_{1,1} = 1$. Then
\[
\Theta(\sigma; \set{1, \ldots, n}, \vec{u}) = c^{k-1} \prod_{i=1}^k \ip{e_{u(b_{i,1})}}{T_{u(b_{i,2})} \ldots T_{u(b_{i, j(i) - 1})} e_{u(b_{i, j(i)})}}.
\]
\end{Lemma}

\begin{proof}
This follows from the definition of $\beta$, noting that $W(b_{1,1}) = W(1) = a_{u(1)}^-$, $W(b_{1,j}) = \tilde{a}_{u(b_{1,l})} = c a_{u(b_{1,l})}^-$ for $j \neq 1$, $W(b_{i, j(l)}) = a_{u(b_{i, j(l)})}^+$, and the rest of the terms are $T_{u(b_{i,l})}$.
\end{proof}

\section{Main theorems}
\label{Section:Meixner}
\begin{Thm}
\label{Thm:Cumulants}
For each $i$, let
\[
S_i = a_i^+ + T_i + \tilde{a}_i = X_i - a_i^-
\]
be an operator on $\Falg(\mc{H})$. Then the free cumulants of the Fock state $\phi_{C, \set{T_i}}$ from Definition~\ref{Defn:Fock-state} are given by the formula $\Cum{x_i} = 0$,
\[
\Cum{x_i P(\mb{x}) x_j} = \ip{e_i}{P(\mb{S}) e_j} = \ip{e_i}{P(\mb{S}) e_j}_C.
\]
\end{Thm}

\begin{proof}
Since $S_{u(i)} = a_{u(i)}^+ + T_{u(i)} + \tilde{a}_{u(i)}$, and using Notation~\ref{Notation:Covered-bijection}, for $\abs{\vec{u}} = n$,
\[
\begin{split}
\ip{e_{u(1)}}{S_{u(2)} \ldots S_{u(n-1)} e_{u(n)}}
& = \ip{\Omega}{a_{u(1)}^- S_{u(2)} \ldots S_{u(n-1)} a_{u(n)}^+ \Omega}
= \sum_{W \in \mc{W}_n'(\vec{u})} \ip{\Omega}{W \Omega} \\
& = \sum_{\sigma \in \NC_0'(n)}\ip{\Omega}{\beta_{\vec{u}}^{-1}(\set{1, \ldots, n}, \sigma) \Omega} \\
& = \sum_{\sigma \in \NC_0'(n)} \Theta(\sigma; \set{1, \ldots, n}, \vec{u}).
\end{split}
\]
Similarly, since $X_{u(i)} = a_{u(i)}^+ + T_{u(i)} + a_{u(i)}^- + \tilde{a}_{u(i)}$, using Lemma~\ref{Lemma:Bijection} and the preceding equation,
\[
\begin{split}
\ip{\Omega}{X_{u(1)} X_{u(2)} \ldots X_{u(n)} \Omega}
& = \sum_{W \in \mc{W}_n(\vec{u})} \ip{\Omega}{W(1) W(2) \ldots W(n) \Omega} \\
& = \sum_{k=1}^n \sum_{\substack{\pi \in \NC_0(n) \\ \pi = (V_1, V_2, \ldots, V_k)}} \sum_{\substack{\sigma_j \in \NC_0'(V_j) \\ j = 1, \ldots, k}} \prod_{i=1}^k \Theta(\sigma_i; V_i, (\vec{u}:V_i)) \\
& = \sum_{k=1}^n \sum_{\substack{\pi \in \NC_0(n) \\ \pi = (V_1, V_2, \ldots, V_k)}} \prod_{i=1}^k \left( \sum_{\sigma_i \in \NC_0'(V_i)} \Theta(\sigma_i; V_i, (\vec{u}:V_i)) \right) \\
& = \sum_{k=1}^n \sum_{\substack{\pi \in \NC_0(n) \\ \pi = (V_1, V_2, \ldots, V_k)}} \prod_{i=1}^k \ip{e_{(\vec{u}:V_i)(1)}}{S_{(\vec{u}:V_i)(2)} \ldots S_{(\vec{u}:V_i)(n-1)} e_{(\vec{u}:V_i)(n)}}
\end{split}
\]
Thus
\[
\state{x_{\vec{u}}} = \sum_{\pi \in \NC_0(n)} \prod_{V \in \pi} \ip{e_{(\vec{u}:V)(1)}}{S_{(\vec{u}:V)(2)} \ldots S_{(\vec{u}:V)(n-1)} e_{(\vec{u}:V)(n)}}.
\]
Since $\Cum{x_i} = \state{x_i} = 0$, the conclusion of the theorem now follows from the defining relation for the free cumulants, namely
\begin{equation*}
\state{x_{\vec{u}}} = \sum_{\pi \in \NC(n)} \prod_{B \in \pi} \Cum{\prod_{i \in B} x_{u(i)}}.
\qedhere
\end{equation*}
\end{proof}

\begin{Cor}
\label{Cor:Self-adjoint}
Each $X_i$ is symmetric and bounded, hence self-adjoint.
\end{Cor}

\begin{proof}
The symmetry is proved exactly as in Proposition~1 of \cite{AnsMonic}, or can be deduced from it. To prove boundedness, choose $m$ such that $\norm{C}, \norm{T_i} < m$. Since $\abs{\NC(n)} < 4^n$, and $\abs{\Theta(\pi; V, \vec{u})} < m^{\abs{V}}$, it follows that $\Cum{x_{\vec{u}}} < (4m)^{\abs{u}}$ and
\[
\phi_{C, \set{T_i}} \left[ X_{u(1)} X_{u(2)} \ldots X_{u(n)} \right] < (16 m)^n.
\]
Thus for each $i$, $\norm{X_i} < 16 m$.
\end{proof}

\begin{Notation}
Let $\mb{z} = (z_1, \ldots, z_d)$ be non-commuting indeterminates, which commute with $\mb{x}$. For a non-commutative power series $G$ in $\mb{z}$ and $i = 1, \ldots, d$, define the left non-commutative partial derivative $D_i G$ by a linear extension of $D_i(1) = 0$,
\[
D_i z_{\vec{u}} = \delta_{i u(1)} z_{u(2)} \ldots z_{u(n)}.
\]
Denote by $\mb{D} G = (D_1 G, \ldots, D_d G)$ the left non-commutative gradient.

\br
For a non-commutative power series $G$, denote by $G^{-1}$ its inverse with respect to multiplication. For a $d$-tuple of non-commutative power series $\mb{G} = (G_1, \ldots, G_d)$, denote by $\mb{G}^{\langle -1 \rangle}$ its inverse with respect to composition (which is also a $d$-tuple).
\end{Notation}

\begin{Thm}
\label{Thm:Meixner}
Let $\phi$ be a state on $\mf{R} \langle \mb{x} \rangle$ with a monic orthogonal polynomial system (MOPS), zero means and identity covariance. The following are equivalent.
\begin{enumerate}
\item
There exists a non-commutative power series
\[
F(\mb{z}) = 1 + (\textsl{terms of degree } \geq 2)
\]
and a $d$-tuple of non-commutative power series $\mb{U}$,
\[
U_i(\mb{z}) = z_i + \textsl{higher-order terms},
\]
such that the polynomials defined via their generating function
\[
\sum_{\abs{\vec{u}} \geq 0} P_{\vec{u}}(\mb{x}) z_{\vec{u}} = F(\mb{z}) \Bigl(1 - \mb{x} \cdot \mb{U}(\mb{z})\Bigr)^{-1}
\]
are a MOPS for $\phi$.
\item
The polynomials with the generating function
\begin{equation}
\label{Generating}
\sum_{\abs{\vec{u}} \geq 0} P_{\vec{u}}(\mb{x}) z_{\vec{u}} = \Bigl( 1 - \mb{x} \cdot (\mb{D} R)^{\langle -1 \rangle} (\mb{z}) + R\bigl((\mb{D} R)^{\langle -1 \rangle} (\mb{z})\bigr) \Bigr)^{-1}
\end{equation}
are a MOPS for $\phi$, where $R$ is the free cumulant generating function~\eqref{Non-crossing} of $\phi$.
\item
The free cumulant generating function of $\phi$ satisfies, for each $i, j$, a (non-commutative) second-order partial differential equation
\begin{equation}
\label{PDE}
D_i D_j R(\mb{z}) = \delta_{ij} + \sum_{k=1}^d B_{ij}^k D_k R(\mb{z}) + C_{ij} D_i R(\mb{z}) D_j R(\mb{z}),
\end{equation}
where $C_{ij} \geq -1$, $B_{ij}^{k} = B_{ik}^{j}$, and for each $j,k$, either $B_{ij}^{k} = 0$ for all $i$, or $C_{ju} = C_{ku}$ for all $u$.
\item
There is a family of polynomials $\set{P_{\vec{u}}}$ such that $\state{P_{\vec{u}}} = 0$ for all $\vec{u} \neq \emptyset$ and they satisfy a recursion relation
\begin{align*}
x_i & = P_i, \\
x_i P_{j} &= P_{(i,j)} + \sum_{k=1}^d B_{ij}^{k} P_{k} + \delta_{ij}, \\
x_i P_{(j, \vec{u})} &= P_{(i, j, \vec{u})} + \sum_{k=1}^d B_{ij}^{k} P_{(k, \vec{u})} + \delta_{ij} (1 + C_{i, u(1)}) P_{\vec{u}},
\end{align*}
where $C_{ij}, B_{ij}^{k}$ satisfy the same conditions as in part (c).
\item
There exist symmetric matrices $T_i$ and a diagonal non-negative matrix $C$ with $(T_i \otimes I) C = C (T_i \otimes I)$ such that $\phi$ has a representation $\phi_{C, \set{T_i}}$ as a Fock state of Definition~\ref{Defn:Fock-state}.
\end{enumerate}
We call such states \emph{free Meixner states}.
\end{Thm}

\begin{proof}
The equivalence (a)~$\Leftrightarrow$~(b) follows from Lemma~4 of \cite{AnsMulti-Sheffer} and Theorem 3.21 of \cite{AnsAppell}, neither of which relied on the assumption that $\phi$ is faithful. The equivalence (d)~$\Leftrightarrow$~(e) follows from the equivalence between the more general Fock space construction and the more general recursion relation in Theorem~\ref{Thm:Monic-states}.

\br
(e) $\Rightarrow$ (c).
By Theorem~\ref{Thm:Cumulants},
\[
R(\mb{z}) = \sum_{j, l = 1}^d \biggl( \ip{e_j}{e_l} z_j z_l + \sum_{\abs{\vec{u}} \geq 1} \ip{e_j}{S_{\vec{u}} e_l} z_j z_{\vec{u}} z_l \biggr).
\]
Therefore
\[
D_j R(\mb{z}) = \sum_{l = 1}^d \biggl( \ip{e_j}{e_l} z_l + \sum_{\abs{\vec{u}} \geq 1} \ip{e_j}{S_{\vec{u}} e_l} z_{\vec{u}} z_l \biggr)
\]
and
\[
\begin{split}
D_i D_j R(\mb{z})
& = \ip{e_j}{e_i} + \sum_{l = 1}^d \biggl( \ip{e_j}{S_i e_l} z_l + \sum_{\abs{\vec{u}} \geq 1} \ip{e_j}{S_i S_{\vec{u}} e_l} z_{\vec{u}} z_l \biggr) \\
& = \ip{e_j}{e_i} + \sum_{l = 1}^d \biggl( \ip{e_j}{T_i e_l} z_l + \sum_{\abs{\vec{u}} \geq 1} \ip{e_j}{(T_i + \tilde{a}_i) S_{\vec{u}} e_l} z_{\vec{u}} z_l \biggr) \\
& = \ip{e_j}{e_i} + \sum_{l = 1}^d \biggl( \ip{e_j}{T_i e_l} z_l + \sum_{\abs{\vec{u}} \geq 1} \ip{e_j}{T_i S_{\vec{u}} e_l} z_{\vec{u}} z_l \biggr)
+ \sum_{l = 1}^d \sum_{\abs{\vec{u}} \geq 1} \ip{e_j}{\tilde{a}_i S_{\vec{u}} e_l} z_{\vec{u}} z_l \\
& = \ip{e_j}{e_i} + \sum_{l = 1}^d \biggl( \sum_{k=1}^d \ip{e_j}{T_i e_k} \ip{e_k}{e_l} z_l + \sum_{\abs{\vec{u}} \geq 1} \sum_{k=1}^d \ip{e_j}{T_i e_k} \ip{e_k}{S_{\vec{u}} e_l} z_{\vec{u}} z_l \biggr) \\
&\quad + \sum_{l = 1}^d \sum_{\abs{\vec{u}} \geq 1} \ip{e_j}{\tilde{a}_i S_{\vec{u}} e_l} z_{\vec{u}} z_l
\end{split}
\]
where in the last step we have used the fact that $\set{e_k}$ form an orthonormal basis. Using Lemma~\ref{Lemma:Factor2}, for $n \geq 4$ and $\vec{u}$ a multi-index on $\set{3, \ldots, n-1}$
\[
\begin{split}
\ip{e_j}{\tilde{a}_i S_{\vec{u}} e_l}
& = \ip{\Omega}{a_j^- \tilde{a}_i S_{\vec{u}} e_l} \\
& = C_{ij} \sum_{k=3}^{n-1}
\sum_{\begin{subarray}{l}
W_1 \in \mc{W}_{\set{1, k+1, \ldots, n}}'\bigl(j, (\vec{u}:\set{k+1, \ldots, n-1}), l\bigr) \\
W_2 \in \mc{W}_{\set{2, \ldots, k}}'\bigl(i, (\vec{u}:\set{3, \ldots, k})\bigr)
\end{subarray}}
\ip{\Omega}{W_1 \Omega} \ip{\Omega}{W_2 \Omega} \\
& = C_{ij} \sum_{k=3}^{n-1} \sum_{W_1 \in \mc{W}_{\set{1, k+1, \ldots, n}}'\bigl((j, \vec{w}, l)\bigr)}
\sum_{W_2 \in \mc{W}_{\set{2, \ldots, k}}'\bigl((i, \vec{v})\bigr)} \ip{e_j}{W_1 e_l} \ip{e_i}{W_2 e_s},
\end{split}
\]
where $\vec{v} = (\vec{u}:\set{3, \ldots, k-1})\bigr)$, $s = u(k)$, and $\vec{w} = (\vec{u}:\set{k+1, \ldots, n})\bigr)$. Thus
\[
\begin{split}
D_i D_j R(\mb{z})
& = \ip{e_j}{e_i} + \sum_{l = 1}^d \left( \sum_{k=1}^d \ip{e_j}{T_i e_k} \ip{e_k}{e_l} z_l + \sum_{\abs{\vec{u}} \geq 1} \sum_{k=1}^d \ip{e_j}{T_i e_k} \ip{e_k}{S_{\vec{u}} e_l} z_{\vec{u}} z_l \right) \\
&\quad + C_{ij} \sum_{l = 1}^d \sum_{(\vec{v}, s, \vec{w})} \sum_{W_1 \in \mc{W}_{\set{1, k+1, \ldots, n}}'\bigl((j, \vec{w}, l)\bigr)}
\sum_{W_2 \in \mc{W}_{\set{2, \ldots, k}}'\bigl((i, \vec{v})\bigr)} \ip{e_j}{W_1 e_l} \ip{e_i}{W_2 e_s} z_{\vec{v}} z_s z_{\vec{w}} z_l \\
& = \ip{e_j}{e_i} + \sum_{k=1}^d \ip{e_j}{T_i e_k} D_k R(\mb{z})
+ \sum_{l = 1}^d  \sum_{\substack{\vec{u} = (\vec{v}, s, \vec{w}) \\
\abs{\vec{v}}, \abs{\vec{w}} \geq 0}} C_{ij} \ip{e_i}{S_{\vec{v}} e_s} \ip{e_j}{S_{\vec{w}} e_l} z_{\vec{v}} z_s z_{\vec{w}} z_l \\
& = \ip{e_j}{e_i} + \sum_{k=1}^d \ip{e_j}{T_i e_k} D_k R(\mb{z}) + C_{ij} D_i R(\mb{z}) D_j R(\mb{z}).
\end{split}
\]
The conditions on the coefficients in part (c) are equivalent to the conditions on the matrices in part (e).

\br
(c) $\Rightarrow$ (e).
Since the states are assumed to have zero means, the corresponding free cumulant generating functions have no linear terms. In that case, a free cumulant generating function $R$, and so the corresponding state $\phi$,
are completely determined by equations~\eqref{PDE}. Moreover, for any choice of $\set{C_{ij}, B_{ij}^k}$ subject to the conditions of part (c), if
\[
T_i (e_j) = \sum_{k=1}^d B_{ij}^k e_k
\]
and
\[
C(e_i \otimes e_j) = C_{ij} \ e_i \otimes e_j,
\]
then those equations are satisfied by $R_{\phi_{C, \set{T_i}}}$. So the states whose free cumulant generating functions satisfy the equations in part (c) are exactly the states in part (e).

\br
(b) $\Rightarrow $ (c).
$\phi$ has a MOPS, so by Theorem~\ref{Thm:Monic-states}, $\phi = \phi_{\mc{C}, \set{\mc{T}_i}}$ for some $\set{\mc{C}^{(k)}, \mc{T}^{(k)}_i}$. Thus, $\phi$ is the joint distribution of the operators $(\mc{X}_1, \ldots, \mc{X}_d)$ on the Hilbert space $\mc{F}_{\mc{C}}(\mc{H})$, with
\[
\mc{X}_i = a_i^+ + \mc{T}_i + a_i^- \mc{C}.
\]
Note that since $\phi$ has means zero and identity covariance, $\mc{T}_i^{(0)} = 0$ and $\mc{C}^{(1)} = I$. Using notation from Section~\ref{Subsubsec:General-Fock}, and denoting
\[
(DR)_{\vec{u}}(\mb{z}) = D_{u(1)} R (\mb{z}) \ldots D_{u(\abs{\vec{u}})} R (\mb{z})
\]
and
\[
e_{\vec{u}} = e_{u(1)} \otimes \ldots \otimes e_{u(\abs{\vec{u}})},
\]
we see that
\[
\begin{split}
& \Bigl(1 - \mb{X} \cdot \mb{z} + R(\mb{z}) \Bigr) \Bigl( \Omega + \sum_{\vec{u}} (DR)_{\vec{u}}(\mb{z}) e_{\vec{u}} \Bigr) \\
&\quad = \Omega - \sum_{i=1}^d z_i e_i + R(\mb{z}) \Omega + \sum_{\vec{u}} (DR)_{\vec{u}}(\mb{z}) e_{\vec{u}} + R(\mb{z}) \sum_{\vec{u}} (DR)_{\vec{u}}(\mb{z}) e_{\vec{u}} \\
&\qquad - \sum_{i=1}^d \sum_{\vec{u}} z_i (DR)_{\vec{u}}(\mb{z}) e_{(i,\vec{u})}
- \sum_{i=1}^d \Bigl(z_i D_i R(\mb{z}) \Omega + \sum_{\vec{u}} z_i (D_i R)(\mb{z}) (DR)_{\vec{u}}(\mb{z}) e_{\vec{u}} \Bigr) \\
&\qquad - \sum_{i}^d \sum_{\vec{u}} z_i (DR)_{\vec{u}}(\mb{z}) \mc{T}_i (e_{\vec{u}})
- \sum_{i}^d \sum_{\vec{u}} z_i (DR)_{\vec{u}}(\mb{z}) a_i^- (\mc{C} - I) e_{\vec{u}}.
\end{split}
\]
Since for any function $G$ with zero constant term,
\begin{equation}
\label{Integral}
\sum_{i=1}^d z_i D_i G(\mb{z}) = G(\mb{z}),
\end{equation}
the preceding expression equals
\[
\begin{split}
& = \Omega - \sum_{i=1}^d z_i e_i + \sum_{\vec{u}} (DR)_{\vec{u}}(\mb{z}) e_{\vec{u}}
- \sum_{i=1}^d \sum_{\vec{u}} z_i (DR)_{\vec{u}}(\mb{z}) e_{(i,\vec{u})} \\
&\quad - \sum_{i}^d \sum_{\vec{u}} z_i (DR)_{\vec{u}}(\mb{z}) \mc{T}_i (e_{\vec{u}})
- \sum_{i}^d \sum_{\vec{u}} z_i (DR)_{\vec{u}}(\mb{z}) a_i^- (\mc{C} - I) e_{\vec{u}}.
\end{split}
\]
Using the expansions \eqref{Expansion-T} and \eqref{Expansion-C} from Theorem~\ref{Thm:Monic-states}, we now continue the equation as
\[
\begin{split}
& = \Omega - \sum_{i=1}^d z_i e_i + \sum_{\vec{u}} (DR)_{\vec{u}}(\mb{z}) e_{\vec{u}}
- \sum_{i=1}^d \sum_{\vec{u}} z_i (DR)_{\vec{u}}(\mb{z}) e_{(i,\vec{u})} \\
&\quad - \sum_{i,j,k=1}^d z_i \Bigl(B_{i, k, j} D_k R(\mb{z}) e_j + \sum_{\vec{u}, \vec{w}} B_{i, (k, \vec{u}), (j, \vec{w})} D_k R(\mb{x}) (DR)_{\vec{u}}(\mb{z}) e_{(j,\vec{w})} \Bigr) \\
&\quad - \sum_{i,j=1}^d z_i \Bigl( (C_{(i,j)} - 1) D_i R(\mb{z}) D_j R(\mb{z}) e_j + \sum_{\vec{u}} (C_{(i,j, \vec{u})} - 1) D_i R(\mb{z}) D_j R(\mb{z}) (DR)_{\vec{u}}(\mb{z}) e_{(j, \vec{u})}
\end{split}
\]
which can be re-organized as
\[
\begin{split}
& = \Omega + \sum_{j=1}^d \Bigl[D_j R(\mb{z}) - \sum_{i=1}^d z_i \Bigl(\delta_{ij} + \sum_{k=1}^d B_{i, k, j} D_k R(\mb{z}) + (C_{(i,j)} - 1) D_i R(\mb{z}) D_j R(\mb{z}) \Bigr) \Bigr] e_j \\
&\quad + \sum_{j=1}^d \sum_{\vec{u}} \Bigl[D_j R(\mb{z}) (DR)_{\vec{u}}(\mb{z}) - \sum_{i=1}^d z_i \Bigl(\delta_{ij}  (DR)_{\vec{u}}(\mb{z}) \\
&\qquad + \sum_{k=1}^d \sum_{\vec{w}} B_{i, (k, \vec{u}), (j, \vec{w})} D_k R(\mb{x}) (DR)_{\vec{w}}(\mb{z}) + (C_{(i,j, \vec{u})} - 1) D_i R(\mb{z}) D_j R(\mb{z}) (DR)_{\vec{u}}(\mb{z})\Bigr)\Bigr] e_{(j, \vec{u})}.
\end{split}
\]
Using equation~\eqref{Integral} again, this equals
\begin{equation}
\label{Intermediate}
\begin{split}
& = \Omega + \sum_{i,j=1}^d \Bigl[D_i D_j R(\mb{z}) - \Bigl(\delta_{ij} + \sum_{k=1}^d B_{i, k, j} D_k R(\mb{z}) + (C_{(i,j)} - 1) D_i R(\mb{z}) D_j R(\mb{z}) \Bigr) \Bigr] z_i e_j \\
&\quad + \sum_{i,j=1}^d \sum_{\vec{u}} \Bigl[D_i D_j R(\mb{z}) (DR)_{\vec{u}}(\mb{z}) - \Bigl(\delta_{ij}  (DR)_{\vec{u}}(\mb{z}) \\
&\qquad + \sum_{k=1}^d \sum_{\vec{w}} B_{i, (k, \vec{u}), (j, \vec{w})} D_k R(\mb{x}) (DR)_{\vec{w}}(\mb{z}) + (C_{(i,j, \vec{u})} - 1) D_i R(\mb{z}) D_j R(\mb{z}) (DR)_{\vec{u}}(\mb{z})\Bigr)\Bigr] z_i e_{(j, \vec{u})}.
\end{split}
\end{equation}
If the polynomials $\set{P_{\vec{u}}}$ with the generating function~\eqref{Generating} from part (b) are orthogonal, then
\[
\sum_{\abs{\vec{u}} \geq 0} P_{\vec{u}}(\mb{x})(DR)_{\vec{u}}(\mb{z}) = \Bigl( 1 - \mb{x} \cdot \mb{z} + R(\mb{z}) \Bigr)^{-1},
\]
and
\begin{equation}
\label{Polynomials-vectors}
P_{\vec{u}} (\mb{X}) \Omega = e_{\vec{u}},
\end{equation}
so that
\begin{equation}
\label{Generating-Omega}
\Bigl(1 - \mb{X} \cdot \mb{z} + R(\mb{z}) \Bigr) \Bigl( \Omega + \sum_{\vec{u}} (DR)_{\vec{u}}(\mb{z}) e_{\vec{u}} \Bigr) = \Omega.
\end{equation}
Equating to zero the coefficient of $z_i e_j$ in equation~\eqref{Intermediate}, we get exactly equation~\eqref{PDE} from part (c), with $B_{ij}^k = B_{i,k,j}$ and $C_{uj} = C_{(i,j)} - 1$. The conditions on the coefficients follow from the general conditions in Theorem~\ref{Thm:Monic-states}.

\br
(e) $\Rightarrow$ (b).
If $\phi = \phi_{C, \set{T_i}}$, it follows that in equation~\eqref{Intermediate}, $B_{i, (k, \vec{u}), (j, \vec{w})} = B_{ij}^k \delta_{\vec{u}, \vec{w}} $ and $C_{i, j, \vec{u}} = 1 + C_{ij}$. Then that expression equals to
\[
\begin{split}
= \Omega + \sum_{i,j=1}^d & \Bigl[D_i D_j R(\mb{z}) - \Bigl(\delta_{ij} + \sum_{k=1}^d B_{ij}^k D_k R(\mb{z}) + C_{ij} D_i R(\mb{z}) D_j R(\mb{z}) \Bigr) \Bigr] \\
& \times z_i \Bigl[e_j + \sum_{\vec{u}} (DR)_{\vec{u}}(\mb{z}) e_{(j, \vec{u})} \Bigr] = \Omega
\end{split}
\]
since part (e) $\Rightarrow $ (c). So equation~\eqref{Generating-Omega} holds. As a result, for polynomials with the generating function~\eqref{Generating},
\[
\Bigl(1 + \sum_{\abs{\vec{u}} \geq 0} P_{\vec{u}}(\mb{X}) z_{\vec{u}} \Bigr) \Omega = \Omega + \sum_{\abs{\vec{u}} \geq 0} z_{\vec{u}} e_{\vec{u}}.
\]
Thus equation~\eqref{Polynomials-vectors} holds, and so the polynomials are orthogonal.
\end{proof}

\subsection{Nontrivial covariance and other extensions}
\label{Subsec:Covariance}
In this section we consider a number of constructions and examples that involve free Meixner states with non-trivial covariances. We still assume that they have zero means (for simplicity); if desired, the means $p_1, \ldots, p_d$ can easily be incorporated into the operator model by considering the operators $(X_1 + p_1, \ldots, X_d + p_d)$ instead, and the corresponding combinatorics will involve all non-crossing partitions $\NC(n)$ rather than the non-crossing partitions without singletons $\NC_0(n)$.

\br
On the other hand, Theorem~\ref{Thm:Monic-states} requires that for any state with MOPS, the covariance matrices have to be diagonal. But now we allow
\begin{equation}
\label{Covariance}
\psi \left[ x_i^2 \right] = t_i.
\end{equation}
Note that degenerate variances $\psi \left[ x_i^2 \right] = 0$ are still not permitted.

\subsubsection{Dilations}
Let $\phi$ be a Meixner state, and fix positive numbers $(t_1, t_2, \ldots, t_d)$. Let $\psi$ be the state defined by the $\mf{R}$-linear extension of
\[
\psi \left[ P(x_1, \ldots, x_d) \right] = \state{P(t_1 x_1, \ldots, t_d x_d)}.
\]
Note that equation~\eqref{Covariance} holds. It is easy to see that if $\set{P_{\vec{u}}}$ is a MOPS for $\phi$, then
\[
Q_{\vec{u}}(x_1, \ldots, x_d) = t_{\vec{u}} P_{\vec{u}}(x_1/t_1, \ldots, x_d/t_d)
\]
is a MOPS for $\psi$.

\br
We now briefly state how the results of Theorem~\ref{Thm:Meixner} get modified for $\psi$. The generating function for the MOPS still has the same ``resolvent'' form, and any state with MOPS and such a generating function arises as a dilation of a free Meixner state.
\[
R_\psi \left[ P(x_1, \ldots, x_d) \right] = R_\phi \left[P(t_1 x_1, \ldots, t_d x_d) \right]
\]
which shows how to modify the differential equation satisfied by the free cumulant generating function. Similarly, the MOPS satisfy the recursion relation
\[
x_i Q_{(j, \vec{u})} = Q_{(i, j, \vec{u})} + \sum_{k=1}^d t_i B_{ij}^{k} Q_{(k, \vec{u})} + \delta_{ij} t_i^2 (1 + C_{i, u(1)}) Q_{\vec{u}}.,
\]
where $\set{B_{ij}^k, C_{ij}}$ were the corresponding coefficients for $\phi$. Finally, suppose that $\phi = \phi_{C, \set{T_i}}$, represented as the joint distribution of $(X_1, X_2, \ldots, X_d)$. In the Hilbert space $\mc{H} = \mf{C}^d$ with an orthonormal basis $\set{e_i}$, let $f_i = t_i e_i$. On the Fock space $\mc{F}_C(\mc{H})$, let
\[
a^+_{e_i} := a^+_i, \quad a_{e_i} := a_i^-, \quad T_{e_i} := T_i,
\]
and extend these definitions $\mf{C}$-linearly to $a^+_f$, $a_f$, $T_f$ for any $f \in \mc{H}$. Let
\begin{equation}
\label{Dilated-operator}
X_{f_i} = a^+_{f_i} + T_{f_i} + a_{f_i}^- + a_{f_i}^- C = t_i X_i.
\end{equation}
Then $\psi$ is the joint distribution of $(X_{f_1}, X_{f_2}, \ldots, X_{f_d})$,
\[
\psi \left[ P(\mb{x}) \right] = \ip{\Omega}{P(\mb{X_f}) \Omega}.
\]

\subsubsection{Free convolution semigroups}
Let $\phi$ be a free Meixner state. For $t > 0$, define a linear functional $\phi^{\boxplus t}$ via its free cumulant functional using relation~\eqref{Cumulants-definition}:
\[
R_{\phi^{\boxplus t}} \left[P(\mb{x}) \right] = t R_\phi[P(\mb{x})].
\]
Note that $\phi^{\boxplus t} \left[ x_i^2 \right] = t$. The notation reflects the fact that
\[
\phi^{\boxplus s} \boxplus \phi^{\boxplus t} = \phi^{\boxplus (s+t)},
\]
where $\boxplus$ is the operation of (additive) free convolution; we will not use this property in the paper. Using the methods of Theorem~\ref{Thm:Cumulants}, it is easy to see that $\phi^{\boxplus t}$ is a state (and so positive) if any only if
\[
t + \min_{i,j} C_{ij} \geq 0,
\]
in other words if $t (I \otimes I) + C \geq 0$. In particular, by assumption~\eqref{C-positive}, $\phi^{\boxplus t}$ is always a state for $t \geq 1$; this is typical behavior for free convolution, as indicated by Corollary 14.13 in \cite{Nica-Speicher-book}. $\phi^{\boxplus t}$ is a state for all $t > 0$ if and only if $C \geq 0$; in this case we say that $\phi$ is \emph{freely infinitely divisible}.

\br
Again, $\phi^{\boxplus t}$ has a MOPS and the generating function for the MOPS still has the same ``resolvent'' form.
\[
D_i D_j R_{\phi^{\boxplus t}} = \delta_{ij} t + \sum_{k=1}^d B_{ij}^k D_k R_{\phi^{\boxplus t}} + (C_{ij}/t) D_i R_{\phi^{\boxplus t}} \ D_j R_{\phi^{\boxplus t}}
\]
and
\[
x_i P_{(j, \vec{u})} = P_{(i, j, \vec{u})} + \sum_{k=1}^d B_{ij}^{k} P_{(k, \vec{u})} + \delta_{ij} (t + C_{i, u(1)}) P_{\vec{u}},
\]
where $\set{B_{ij}^k, C_{ij}}$ were the corresponding coefficients for $\phi$. Finally, suppose that $\phi = \phi_{C, \set{T_i}}$. On the algebraic Fock space $\Falg(\mc{H})$, define an inner product using the kernel
\[
\begin{split}
K_C^{(t)} & = \bigl(I^{\otimes (k-2)} \otimes (t I^{\otimes 2} + C)\bigr) \ldots \bigl(I \otimes (t I^{\otimes 2} + C) \otimes I^{\otimes (k-3)}\bigr) \bigl((t I^{\otimes 2} + C) \otimes I^{\otimes (k-2)}\bigr) t \\
& = t^k K_{C/t}
\end{split}
\]
on $\mc{H}^{\otimes k}$, and denote the completion of $\Falg(\mc{H})$ with respect to this inner product $\mc{F}_C^{(t)}(\mc{H})$. Let
\begin{equation}
\label{Xs}
X^{(t)}_i = a^+_{i} + T_{i} + t a_{i}^- + \tilde{a}_i
= a^+_{i} + T_{i} + t a_i^- (I + C/t).
\end{equation}
Then $\phi^{\boxplus t}$ is the joint distribution of $\left(X^{(t)}_1, \ldots, X^{(t)}_d \right)$.

\br
As constructed above, $\set{X^{(t)}_i}$ are represented on different Hilbert spaces for different $t$. We can combine this construction with an idea from Section 7.2 of~\cite{Sniady-SWN} to represent a whole family of functionals $\set{\phi^{\boxplus t} | 0 < t < 1}$ on a single space.

\br
A subset $S \subset \set{1, 2, \ldots, n-1}$ can be identified with an \emph{interval partition} $\pi(S) \in \Int(n)$: if $S = \set{i(1), i(2), \ldots, i(k)}$, then
\[
\pi = \bigl( \set{1, \ldots, i(1)}, \set{i(1) + 1, \ldots, i(2)}, \ldots, \set{i(k) + 1, \ldots, n} \bigr).
\]
Consider the vector space $H = \mc{H} \otimes L^\infty([0,1], dx)$ as a subspace of the Hilbert space $\mc{H} \otimes L^2([0,1], dx)$,
with the inner product
\[
\ip{\eta \otimes f}{\zeta \otimes g} = \ip{\eta}{\zeta} \int_0^1 f(x) g(x) \,dx.
\]
On its algebraic Fock space $\Falg(H)$, define the inner product
\begin{multline*}
\ip{(\eta_1 \otimes f_1) \otimes \ldots \otimes (\eta_l \otimes f_l)}{(\zeta_1 \otimes g_1) \otimes \ldots \otimes (\zeta_n \otimes g_n)}_C \\
= \delta_{ln} \sum_{\substack{S \subset \set{1, \ldots, n-1} \\ \pi(S) = (V_1, V_2, \ldots, V_k)}}
\ip{\eta_1 \otimes \ldots \otimes \eta_n}{C^{S^c} \left(\zeta_1 \otimes \ldots \otimes \zeta_n \right)} \prod_{j=1}^{k} \left( \int_{\mf{R}} \left[ \prod_{i \in V_j} f_i(x) g_i(x) \right] \,dx \right),
\end{multline*}
where $S^c$ is the complement $\set{1, \ldots, n-1} \backslash S$, and
\[
C^{S^c} = \prod_{i \in S^c} I^{\otimes (i-1)} \otimes C \otimes I^{\otimes (n-i-1)}.
\]
Complete with respect to this inner product, to get the Hilbert space $\mc{F}_C (H)$. On this space, define operators
\begin{align*}
{a_i^+}^{(t)} \bigl((\eta_1 \otimes f_1) \otimes \ldots \otimes (\eta_n \otimes f_n)\bigr) & = (e_i \otimes \chf{[0,t)}) \otimes (\eta_1 \otimes f_1) \otimes \ldots \otimes (\eta_n \otimes f_n) \\
{a_i^-}^{(t)} \bigl((\eta_1 \otimes f_1) \otimes \ldots \otimes (\eta_n \otimes f_n)\bigr) & = \ip{e_i}{\eta_1} \Bigl( \int_0^t f_1(x) \,dx \Bigr) (\eta_2 \otimes f_2) \otimes \ldots \otimes (\eta_n \otimes f_n), \\
T_i^{(t)} \bigl((\eta_1 \otimes f_1) \otimes \ldots \otimes (\eta_n \otimes f_n)\bigr) & = (T_i \eta_1 \otimes f_1 \chf{[0,t)}) \otimes (\eta_2 \otimes f_2) \otimes \ldots \otimes (\eta_n \otimes f_n), \\
\tilde{a}_i^{(t)} \bigl((\eta_1 \otimes f_1) \otimes \ldots \otimes (\eta_n \otimes f_n)\bigr) & \\
= \bigl((a_i^- C(\eta_1 & \otimes \eta_2)) \otimes (f_1 \chf{[0,t)} f_2)\bigr) \otimes (\eta_3 \otimes f_3) \otimes \ldots \otimes (\eta_n \otimes f_n),
\end{align*}
where $\chf{[0,t)}$ is the indicator function of the interval $[0,t)$, and let
\[
X_i^{(t)} = {a_i^+}^{{(t)}} + T_i^{(t)} + {a_i^-}^{(t)} + \tilde{a}_i^{(t)}.
\]
By combining Corollary~\ref{Cor:Self-adjoint} with (a slight modification of) Theorem~6 from \cite{Sniady-SWN}, it follows that each $X_i^{(t)}$ is self-adjoint on $\mc{F}_C(H)$. Note that if all $f_i = g_i = \chf{[0,t)}$, then
\[
\begin{split}
& \ip{(\eta_1 \otimes f_1) \otimes \ldots \otimes (\eta_n \otimes f_n)}{(\zeta_1 \otimes g_1) \otimes \ldots \otimes (\zeta_n \otimes g_n)}_C \\
&\qquad = \sum_{\substack{S \subset \set{1, \ldots, n-1} \\ \pi(S) = (V_1, V_2, \ldots, V_k)}}
\ip{\eta_1 \otimes \ldots \otimes \eta_n}{C^{S^c} \left(\zeta_1 \otimes \ldots \otimes \zeta_n \right)} t^k \\
&\qquad = t^n \ip{\eta_1 \otimes \ldots \otimes \eta_n}{\zeta_1 \otimes \ldots \otimes \zeta_n}_{C/t}.
\end{split}
\]
Moreover, each $X_i^{(t)}$ restricted to
\[
\mc{F}_C(\mc{H} \otimes \Span{\chf{[0,t)}}) \cong \mc{F}_C^{(t)}(\mc{H})
\]
is given by the equation~\eqref{Xs}, and so $\phi^{\boxplus t}$ is the joint distribution of $\left(X^{(t)}_1, \ldots, X^{(t)}_d \right)$.

\subsubsection{Rotations}
\label{Subsubsec:Rotations}
Let $O = (O_{ij})$ be an orthogonal $d \times d$ matrix. Let
\[
O^T \mb{x} = \left( \sum_{i=1}^d O_{i1} x_i, \ldots, \sum_{i=1}^d O_{id} x_i \right)
\]
and
\begin{equation}
\label{Change-of-variable}
\phi^O \left[ P(\mb{x}) \right] = \state{P(O^T \mb{x})}.
\end{equation}
We call $\phi^O$ a rotation of $\phi$. $\phi^O$ is the joint distribution of $(X_{f_1}, \ldots, X_{f_d})$ from~\eqref{Dilated-operator}, where we take
\[
f_j = O (e_j) = \sum_{i=1}^d O_{ij} e_i.
\]
$\phi^O$ need not have a MOPS, since the matrix $C$ need not be diagonal in the basis $\set{f_1, \ldots, f_d}$. In fact, it follows from Lemma~9 of \cite{AnsMulti-Sheffer} that $\phi^O$ has a MOPS for \emph{all} $O$ if and only if $C_{ij} = c$ for all $i, j$, and that in this case $\phi^O$ is also a free Meixner state. It is easy to see that more generally, if $S \subset \set{1, \ldots, d}$ and $C_{ij} = c$ for all $i, j \in S$, then $\phi^O$ is a free Meixner state whenever $O (e_k) = e_k$ for all $k \not \in S$.

\subsubsection{Linear transformations}
Finally, one can consider a general invertible change of variables
\[
A^T \mb{x} = \left( \sum_{i=1}^d A_{i1} x_i, \ldots, \sum_{i=1}^d A_{id} x_i \right)
\]
and the corresponding state $\phi^A$ defined as in equation~\eqref{Change-of-variable}. $\phi^A$ is the joint distribution of operators from~\eqref{Dilated-operator}, where we take $f_j = A(e_j) = \sum_{i=1}^d A_{ij} e_i$. As an alternative to our definition, one can call free Meixner states all states obtained by a linear transformation of a free Meixner state with MOPS (compare with \cite{Pommeret-Test}).

\section{Examples}
\label{Section:Examples}

\subsection{Free products}
In preparation for the examples in this section, for the reader's convenience we explain a key notion from free probability. Again, see \cite{VDN,Nica-Speicher-book} for more details.

\br
Let $\phi_1, \ldots, \phi_d$ be one-dimensional states on $\mf{R}[x_1], \ldots, \mf{R}[x_d]$, respectively. There is a canonical way to define their \emph{free product state} $\phi$ on $\mf{R} \langle x_1, \ldots, x_d \rangle$. Combinatorially, a natural way to define $\phi$ is via its MOPS. Let $\set{P_n^{(i)}}$ be the MOPS for $\phi_i$. For a multi-index $\vec{u}$, decompose
\[
x_{\vec{u}} = x_{v(1)}^{i(1)} x_{v(2)}^{i(2)} \ldots x_{v(k)}^{i(k)},
\]
where the consecutive indices $v(j) \neq v(j+1)$, although non-consecutive indices may coincide. Then the MOPS $\set{P_{\vec{u}}}$ for $\phi$ are defined by
\[
P_{\vec{u}}(\mb{x}) = \prod_{j=1}^k P_{i(j)}^{(v(j))}(x_{v(j)}).
\]
For example,
\[
P_{1,1,2,1,2}(\mb{x}) = P_2^{(1)}(x_1) P_1^{(2)}(x_2) P_1^{(1)}(x_1) P_1^{(2)}(x_2).
\]
Note that if one considers polynomials in commuting variables and assumes that \emph{all} $v(j)$ above are different, one gets the usual (Cartesian) product of measures. Also, $\phi$ is a free product state if any only if the elements $x_1, x_2, \ldots, x_d$ are freely independent with respect to $\phi$, in the sense of Voiculescu. This can be taken as the definition of free independence; note that for random variables independent in the usual probabilistic sense, their joint distribution is a product measure. Finally, the crucial property of free cumulant generating functions is their relation to free products: a state $\phi$ is a free product state of $\phi_1, \ldots, \phi_d$ if any only if the free cumulant generating function of $\phi$ decomposes as
\[
R_\phi(\mb{z}) = \sum_{i=1}^d R_{\phi_i}(z_i).
\]
This is often stated as the ``mixed free cumulants are zero'' condition. It is related to the familiar property that the Fourier transform of the joint distribution of independent random variables is the product of their individual Fourier transforms.

\br
It is easy to see that free product free Meixner states are exactly the free products of one-dimensional free Meixner states, see Remark~6 of \cite{AnsMulti-Sheffer}. Recall that these one-dimensional distributions, as described in that remark, Theorem 4 of \cite{AnsMeixner} and Section 2.2 of \cite{Boz-Bryc}, are known. With variance $t$, they are: the semicircular (free Gaussian) distributions $\frac{1}{2 \pi} \sqrt{4 t - x^2} \,dx$, the Marchenko-Pastur (free Poisson) distributions $\frac{1}{2 \pi} \frac{\sqrt{4t - (x-b)^2}}{1 + (b/t) x} \,dx + \text{ possibly one atom}$, and more generally
\[
\frac{1}{2 \pi} \frac{\sqrt{4 (t + c) - (x - b)^2}}{1 + (b/t) x + (c/t^2) x^2} \,dx + \text{ zero, one, or two atoms},
\]
depending on the particular values of $b,c,t$.

\subsection{Semicircular systems}
\label{Subsec:Semicircular}
Let $C = 0$ and all $T_i = 0$. Then
\[
S_i = a_i^+
\]
and
\[
\Cum{x_i x_{\vec{u}} x_j} = \ip{e_i}{S_{\vec{u}} e_j} = 0
\]
for $\abs{\vec{u}} \geq 1$. Thus in distribution, all $S_i \sim 0$. Only second-order free cumulants of $(X_1, X_2, \ldots, X_d)$ are non-zero, and $\phi$ is the distribution of a freely independent semicircular system, the free analog of the standard $d$-dimensional Gaussian distribution.

\subsection{Free Poisson states}
Let $C = 0$ and $T_i$ arbitrary. Then
\[
S_i = a_i^+ + T_i
\]
and
\[
\Cum{x_i x_{\vec{u}} x_j} = \ip{e_i}{S_{\vec{u}} e_j} = \ip{e_i}{T_{\vec{u}} e_j}.
\]
Thus in distribution, $(S_1, S_2, \ldots, S_d) \sim (T_1, T_2, \ldots, T_d)$. It is appropriate to say that in this case, the joint distribution $\phi$ of $(X_1, X_2, \ldots, X_d)$ is $d$-dimensional free Poisson. In \cite{AnsMulti-Sheffer} we showed that if $\phi$ is tracial, then $\phi$ is a rotation of a free product of one-dimensional free Poisson distributions. Whether or not $\phi$ is tracial, the vector $\Omega$ is cyclic and separating for the von Neumann algebra $W^\ast(X_1, X_2, \ldots, X_d)$.

\subsection{Free product states}
\label{Subsec:Free-products}
For two vectors $f, g$, denote by $E_{f,g}$ the corresponding rank one operator,
\[
E_{f,g}(h) = f \ip{g}{h}.
\]
For an orthonormal basis $\set{f_i}$, $E_{f_i, f_j}$ are the corresponding matrix units. For the standard basis $\set{e_i}$, we will denote these simply by $E_{ij}$. In particular, $E_{ii}$ is the orthogonal projection onto $e_i$.

\br
Let $C(e_i \otimes e_j) = c_i \delta_{ij} (e_i \otimes e_j)$, and let $T_i = b_i E_{ii}$. Then $S_i$ acts entirely on the subspace $\mc{F}_{c_i}(\Span{e_i})$, on which it equals
\[
S_i = a_i^+ + b_i + c_i a_i^-.
\]
$a_i^+ + c_i a_i$ has the centered semicircular distribution with variance $c_i$ (note that on $\mc{F}_{C}(\mc{H})$, this operator is not symmetric, so its \emph{star}-distribution is different from the semicircular one). Therefore $S_i$ has the semicircular distribution with mean $b_i$ and variance $c_i$. Also, it follows that
\[
\Cum{x_i x_{\vec{u}} x_j} = \ip{e_i}{S_{\vec{u}} e_j} = 0
\]
unless
\[
i = u(1) = \ldots = u(n) = j.
\]
In other words, all the mixed free cumulants of $(X_1, X_2, \ldots, X_d)$ are zero. This says precisely that their joint distribution $\phi$ is a free product of the distributions of each of $X_1, X_2, \ldots, X_d$. Each of these, in turn, is a one-dimensional free Meixner distribution, whose free cumulants are, up to a shift of index, the moments of the semicircular distribution with mean $b_i$ and variance $c_i$.

\subsection{Exponentiated semicircular systems}
Let $C(e_i \otimes e_j) = c_i (e_i \otimes e_j)$ and $T_i = b_i I$. Then
\[
S_i = a_i^+ + c_i a_i^- + b_i I,
\]
again the distribution of $S_i$ is the semicircular distribution with mean $b_i$ and variance $c_i$, but now the operators $S_i$ themselves are freely independent with respect to the state $\phi$, so that their joint distribution is a free product. The joint distribution of $(X_1, X_2, \ldots, X_d)$ is \emph{not} a free product (typically, not even tracial); it was described in the last section of \cite{AnsMulti-Sheffer}.

\br
Note that the preceding two examples make sense for $-1 \leq c_i < 0$, except that one loses the interpretation of $S_i$ as having a semicircular distribution with variance $c_i$, and the resulting states are not freely infinitely divisible.

\begin{Remark}
The constructions in the previous two examples coincide in the one-dimensional case. That case, and in particular the corresponding free cumulants, were also considered in \cite{AnsMeixner} and described completely in \cite{Boz-Bryc}. Moreover, many one-dimensional free Meixner distributions arise as limits in the central and Poisson limit theorems for the $t$-transformed free convolution, in the sense of \cite{Boz-Wys}; that paper also contains a Fock space construction which coincides with the one-dimensional version of the one in Section~\ref{Subsec:Fock2}.
\end{Remark}

\subsection{Free multinomial states}
\label{Subsec:Free-multinomial}
It is well known that the Bernoulli distribution
\[
(1-p) \delta_0 + p \delta_1
\]
is a Meixner distribution. It was noted in \cite{AnsMeixner} that it is also a (one-dimensional) free Meixner distribution. Moreover, the binomial distributions, which are convolution powers
\[
\bigl((1-p) \delta_0 + p \delta_1\bigr)^{\ast n} = \sum_{k=0}^n \binom{n}{k} (1-p)^{n-k} p^k \delta_k
\]
of the Bernoulli distribution, are all Meixner, and the free binomial distributions, which are free convolution powers $\bigl((1-p) \delta_0 + p \delta_1\bigr)^{\boxplus n}$ of the Bernoulli distribution are free Meixner. In fact, it was noted in \cite{Boz-Bryc} that $\bigl((1-p) \delta_0 + p \delta_1\bigr)^{\boxplus t}$ are free Meixner for all real $t \geq 1$.

\br
It is also well-known that the multinomial distributions are Meixner \cite{Pommeret-Test}. In particular, the basic multinomial distribution
\begin{equation}
\label{Multinomial}
p_1 \delta_{e_1} + p_2 \delta_{e_2} + \ldots p_d \delta_{e_d}
\end{equation}
on $\mf{R}^d$ has this property. We now show that it, and so the free semigroup it generates, also induce free Meixner states. In this example, it is natural to consider the state with non-trivial means and a non-diagonal covariance matrix; an actual free Meixner state can be obtained from it by an affine transformation as in Section~\ref{Subsec:Covariance}.

\br
In the Fock space construction of Section~\ref{Subsec:Fock2}, take $\dim \mc{H} = d-1$ rather than $d$, and
\[
C(e_i \otimes e_j) = - e_i \otimes e_j
\]
for all $i, j$, so that $C_{ij} = -1$. In this case the induced inner product on the Fock space $\Falg(\mc{H})$ is degenerate, and the vector space factors through to simply $\mc{F}_C(\mc{H}) = \mf{C} \oplus \mc{H}$. \emph{In this example only}, choose (linearly dependent) vectors $\set{e_i | i = 1, 2, \ldots d}$ in $\mc{H}$ that are not orthonormal, but instead satisfy
\begin{align*}
\ip{e_i}{e_i} & = p_i (1 - p_i), \\
\ip{e_i}{e_j} & = - p_i p_j,
\end{align*}
where
\[
p_i > 0, \qquad i = 1, 2, \ldots, d, \qquad p_1 + p_2 + \ldots + p_{d} = 1.
\]
Since these numbers are the covariances of the centered version of the basic multinomial distribution \eqref{Multinomial}, the corresponding matrix is positive semi-definite and so the $\set{e_i}$ can be chosen in this fashion.

\br
Let
\begin{align*}
T_i(e_i) & = (1 - 2 p_i) e_i, \\
T_i(e_j) & = - p_i e_j - p_j e_i,
\end{align*}
and define
\[
X_i = a_i^+ + T_i + a_i^- + a_i^- C
\]
as usual, except that $a_i^+ = 0$ on $\mc{H}$. In other words, for $Y_i = X_i + p_i$,
\begin{equation}
\label{Multinomial-operators}
\begin{split}
Y_i \Omega & = e_i + p_i \Omega, \\
Y_i e_i & = (1 - p_i) [e_i + p_i \Omega], \\
Y_i e_j & = - p_j [e_i + p_i \Omega].
\end{split}
\end{equation}

\begin{Prop}
$Y_i$ is an orthogonal projection of $\mf{C} \oplus \mc{H}$ onto $\Span{e_i + p_i \Omega}$. These projections are orthogonal among themselves and their sum is the identity operator. Their joint distribution with respect to the state $\phi_{C, \set{T_i}}$ from Definition~\ref{Defn:Fock-state} is the basic multinomial distribution \eqref{Multinomial}. In particular, $\state{Y_i} = p_i$. The free cumulant generating function of $\phi$ satisfies the differential equation
\[
\begin{split}
D_i D_j R
& = (\delta_{ij} p_i - p_i p_j) + (\delta_{ij} - p_j) D_i R - p_i D_j R - D_i R \ D_j R \\
& = \delta_{ij} (D_i R + p_i) - (D_i R + p_i) (D_j R + p_j).
\end{split}
\]
\end{Prop}

\begin{proof}
$Y_i$ is self-adjoint, its image is $\Span{e_i + p_i \Omega}$, and
\[
Y_i(e_i + p_i \Omega)
= (1 - p_i) [e_i + p_i \Omega] + p_i [e_i + p_i \Omega]
= e_i + p_i \Omega,
\]
so it is an orthogonal projection onto its image.
\[
\ip{e_i + p_i \Omega}{e_j + p_i \Omega}
= - p_i p_j + p_i p_j = 0,
\]
so these subspaces, and therefore projections onto them, are orthogonal. $\mf{C} \oplus \mc{H}$ has dimension $d$, therefore the sum $\sum_{i=1}^{d} Y_i$ is identity. It also follows that
\[
\state{Y_{\vec{u}}} =
\begin{cases}
p_i, & i = u(1) = u(2) = \ldots, \\
0, & \text{ otherwise}.
\end{cases}
\]
Therefore, the joint distribution of $(Y_1, Y_2, \ldots, Y_{d})$ with respect to $\phi$ is the basic multinomial distribution. The last part follows from the operator representation.
\end{proof}

\begin{Defn}
\label{Defn:Free-multinomial}
Free multinomial states are the states $\set{\phi^{\boxplus t} | t \geq 1}$, where $\phi$ is the basic multinomial distribution~\eqref{Multinomial}. Note that $\phi^{\boxplus n}$ is the joint distribution of the sum of $n$ $d$-tuples of orthogonal projections.
\end{Defn}

\begin{Remark}
Since the Bernoulli distribution is both classical and free Meixner, one may conjecture that it is in some sense also $q$-Meixner (for $0 \leq q \leq 1$, with $q=1$ corresponding to the classical case and $q=0$ corresponding to the free case). The meaning of this term is not well-defined, but see for example Section 4.3 of \cite{AnsAppell}. Indeed, the recursion relation for its orthogonal polynomials is of the $q$-Meixner form
\begin{align*}
x P_0 & = P_1 + p, \\
x P_1 & = P_2 + (1-p) P_1 + p (1-p) P_0, \\
x P_n & = P_{n+1} + (1-p) [n]_q P_n + [n]_q p (1-p)(1 - [n-1]_q) P_{n-1}
\end{align*}
independently of $q$, as long as the degree of the polynomial $n \leq 1$, which suffices since
\[
L^2\bigl((1-p)\delta_0 + p \delta_1\bigr)
\]
is $2$-dimensional. One may also hope that its $q$-cumulant generating function would then satisfy the equation
\[
D_q^2 R^{(q)} = D_q R^{(q)} - (D_q R^{(q)})^2,
\]
where $D_q$ is the $q$-derivative
\[
D_q(f)(z) = \frac{f(z) - f(qz)}{(1-q) z},
\]
and
\[
R^{(q)} = \sum_{n=1}^\infty \frac{1}{[n]_q!} \alpha_n z^n.
\]
The corresponding recursion for its $q$-cumulants is
\[
\alpha_{n+2} = \alpha_{n+1} - \sum_{i=0}^n \left[ \begin{matrix} n \\ i\end{matrix} \right]_q \alpha_{i+1} \alpha_{n-i+1}
\]
(compare with Remark 5.4 of \cite{Boz-Bryc}), with the initial condition $\alpha_1 = p$. Using Maple, it is easy to calculate the first $5$ cumulants. Unfortunately, the fifth $q$-cumulant of the Bernoulli distribution calculated in this fashion differs from its fifth $q$-cumulant in the sense of Section 6 of \cite{AnsQCum}.
\end{Remark}

\subsection{Tracial examples}
If $\phi$ is a state on a non-commutative algebra $\mc{A}$, one says that $\phi$ is \emph{tracial}, or \emph{a trace}, if for any $x, y \in \mc{A}$,
\[
\state{x y} = \state{y x}.
\]
Tracial states play a crucial role, for example, in the theory of von Neumann algebras.

\begin{Lemma}
Let $\phi = \phi_{C, \set{T_i}}$ be a free Meixner state, represented as the joint distribution of operators $(X_1, \ldots, X_d)$. Suppose that $\phi$ is tracial. Then for all $i, j$,
\begin{equation}
\label{Trace1}
T_i e_j = T_j e_i
\end{equation}
and
\begin{equation}
\label{Commutator}
T_i T_j - T_j T_i = C_{ji} E_{ij} - C_{ij} E_{ji}.
\end{equation}
\end{Lemma}

\begin{proof}
Since $\phi$ is tracial, for all $i, j, k$,
\[
\state{X_i X_j X_k} = \ip{e_i}{T_j e_k} = \ip{e_j}{T_k e_i} = \ip{e_i}{T_k e_j},
\]
so for all $j, k$,
\[
T_j e_k = T_k e_j.
\]
Similarly, for all $i, j, k, l$,
\[
\begin{split}
\state{X_i X_j X_k X_l}
& = \ip{e_i}{e_j} \ip{e_k}{e_l} + \ip{e_i}{e_l} \ip{e_j}{e_k} (1 + C_{kl}) + \ip{e_i}{T_j T_k e_l} \\
& = \ip{e_j}{e_k} \ip{e_l}{e_i} + \ip{e_j}{e_i} \ip{e_k}{e_l} (1 + C_{li}) + \ip{e_j}{T_k T_l e_i},
\end{split}
\]
so
\[
\ip{e_i}{e_l} \ip{e_j}{e_k} C_{kl} + \ip{e_i}{T_j T_k e_l}
= \ip{e_j}{e_i} \ip{e_k}{e_l} C_{li} + \ip{e_i}{T_l T_k e_j}.
\]
Using equation~\eqref{Trace1} and the orthonormality of $\set{e_i}$,
\[
\ip{e_i}{e_l} \ip{e_j}{e_k} C_{jl} + \ip{e_i}{T_j T_l e_k}
= \ip{e_j}{e_i} \ip{e_k}{e_l} C_{lj} + \ip{e_i}{T_l T_j e_k},
\]
so
\[
\ip{e_j}{e_k} C_{jl} e_l + T_j T_l e_k
= \ip{e_k}{e_l} C_{lj} e_j + T_l T_j e_k
\]
and
\[
T_j T_l - T_l T_j = C_{lj} E_{jl} - C_{jl} E_{lj}.
\qedhere
\]
\end{proof}

\begin{Ex}
A general theorem of Voiculescu (Proposition~2.5.3 in \cite{VDN}) implies that the free product states from Example~\ref{Subsec:Free-products} are tracial. It is also easy to see that any rotation of a tracial state is tracial.
\end{Ex}

\begin{Lemma}
\label{Lemma:Tracial-semigroup}
If $\phi$ is tracial, then $\phi^{\boxplus t}$ is tracial for all $t$ for which it is defined.
\end{Lemma}

\begin{proof}
It is easy to see directly from the defining equation~\eqref{Cumulants-definition} that a state $\phi$ is tracial if any only if its free cumulant generating functional $R_\phi$ is tracial. The combinatorial reason is that if a partition $\pi$ is non-crossing when points $\set{1, 2, \ldots, n}$ are placed on a line, it is also non-crossing when they are placed on a circle. The lemma follows from this fact and the defining relation for $\phi^{\boxplus t}$
\[
R_{\phi^{\boxplus t}} \left[ x_{\vec{u}} \right] = t R_{\phi} \left[ x_{\vec{u}} \right].
\qedhere
\]
\end{proof}

\begin{Prop}
All free multinomial states are tracial.
\end{Prop}

\begin{proof}
Since the operators $\set{Y_i}$ defined in equation~\eqref{Multinomial-operators} commute, the basic multinomial distribution is tracial; in fact, it factors through to a state on commutative polynomials $\mf{R}[x_1, \ldots, x_d]$ corresponding to the basic multinomial measure~\eqref{Multinomial}. It follows from Lemma~\ref{Lemma:Tracial-semigroup} that all the other free multinomial states are tracial as well.
\end{proof}

\noindent
We conclude the paper with three further results on when Meixner states are traces. The first two show that under one set of general assumptions, the only tracial Meixner states are the trivial ones, namely the rotations of free product states. It generalizes Proposition~11 of \cite{AnsMulti-Sheffer}. The last one provides a way to construct a large class of tracial examples that do not come from free products. It generalizes the multinomial example above.

\begin{Prop}
\label{Prop:Free-products}
Let $\phi = \phi_{C, \set{T_i}}$ be a tracial free Meixner state with $C$ diagonal as a $d \times d$ matrix, $C_{ij} = \delta_{ij} c_i$. Then $\phi$ is a rotation of a free product state.
\end{Prop}

\begin{proof}
If $C_{ij} = \delta_{ij} c_i$, then equation~\eqref{Commutator} states that all $T_i, T_j$ commute. Combining this with equation~\eqref{Trace1}, we see that moreover, for some orthonormal basis $\set{f_1, f_2, \ldots, f_d}$,
\begin{equation}
\label{Diagonalized}
T_j = \sum_{i=1}^d \alpha_i \ip{f_i}{e_j} E_{f_i, f_i}.
\end{equation}
From
\[
C (T_j \otimes I) (e_k \otimes e_l)
= (T_j \otimes I) C (e_k \otimes e_l)
\]
it follows that
\[
\sum_{i=1}^d \alpha_i \ip{f_i}{e_j} \ip{f_i}{e_k} \ip{f_i}{e_l} c_l (e_l \otimes e_l)
= \delta_{kl} \sum_{i,m=1}^d \alpha_i \ip{f_i}{e_j} \ip{f_i}{e_k} \ip{f_i}{e_m} c_l (e_m \otimes e_l).
\]
Thus for all $k \neq l$,
\[
\sum_{i=1}^d \alpha_i \ip{f_i}{e_j} \ip{f_i}{e_k} \ip{f_i}{e_l} c_l
= \ip{e_k}{\left( \sum_{i=1}^d \alpha_i \ip{f_i}{e_j} E_{f_i, f_i} \right) e_l} c_l
= \ip{e_k}{T_j e_l} c_l = 0.
\]
It follows that whenever $c_l \neq 0$, $T_j e_l \in \Span{e_l}$, and one can take $f_l = e_l$. So if $S = \set{l | c_l \neq 0}$, then
\[
\Span{e_l | l \in S} = \Span{f_l | l \in S}
\]
is an invariant subspace for all $T_j$. We can choose an orthogonal transformation $O$ so that $O (e_i) = f_i$, in other words
\[
O (e_i) =
\begin{cases}
e_i & \text{ for } i \in S, \text{ that is } c_i \neq 0, \\
f_i & \text{ for } i \not \in S, \text{ that is } c_i = 0.
\end{cases}
\]
Following the comments at the end of Section~\ref{Subsubsec:Rotations}, the state $\phi^O$, which is the joint distribution of $(X_{f_1}, \ldots, X_{f_d})$, is still a tracial free Meixner state. From equation~\eqref{Diagonalized},
\[
T_{f_j} = \alpha_j E_{f_j, f_j}.
\]
Finally, $C(\eta \otimes \zeta) = 0$ whenever one of $\eta, \zeta \in \Span{e_l | l \not \in S} = \Span{f_l | l \not \in S}$, so
\[
C(f_i \otimes f_j) =
\begin{cases}
C(e_i \otimes e_j) = c_i \delta_{ij} & \text{ if } i, j \in S, \\
0 & \text{ if one of } i, j \not \in S.
\end{cases}
\]
Thus $C, \set{T_i}$ have the form in Example~\ref{Subsec:Free-products}, and so $\phi^O$ is a free product state.
\end{proof}

\begin{Cor}
Let $\phi = \phi_{C, \set{T_i}}$ be a tracial free Meixner state and $T_i = 0$ for all $i$. Then $\phi$ is a free product state.
\end{Cor}

\begin{proof}
If $T_i = 0$ for all $i$, then it follows from equation~\eqref{Commutator} that $C_{ij} = \delta_{ij} c_i$, and the preceding proposition applies. In this case, a rotation is unnecessary.
\end{proof}

\noindent
Meixner states correspond to quadratic natural exponential families. The following states are free versions of \emph{simple} quadratic natural exponential families in the terminology of \cite{Casalis-Simple-quadratic}, where all such (classical) families were classified.

\begin{Prop}
Let $C$ be a constant matrix, $C_{ij} = c$ for all $i,j$. Then the necessary conditions \eqref{Trace1} and~\eqref{Commutator} for $\phi$ to be a tracial free Meixner state, namely that $\set{T_i}$ are symmetric matrices, $T_i e_j = T_j e_i$ and
\[
(T_i T_j - T_j T_i)  = c \bigl( E_{ij} - E_{ji} \bigr),
\]
are also sufficient.
\end{Prop}

\begin{proof}
If $c=0$, the result follows from Proposition~\ref{Prop:Free-products}. So we will assume that $c \neq 0$.

\br
For each $n$, let
\[
A(n) = \set{\pi \in \NC_0'(n) | 1 \stackrel{\pi}{\sim} 2}.
\]
For any partition $\sigma \in \NC_0'(n) \backslash A(n)$, define the partition $l(\sigma) \in A(n)$ as follows: if $1 \in B \in \sigma$ and $2 \in C \in \sigma$, let $l(\sigma)$ be the partition with the same classes as $\sigma$ except that $B \cup C$ is a class of $l(\sigma)$. Conversely, any such $\sigma$ can be obtained by starting with $\pi \in A(n)$ with the (unique) outer class $B$, choosing $i \in B$, $2 < i < n$ (if it exists), and taking $\sigma$ to have the same classes as $\pi$ except that
\[
B \cap \left( \set{1} \cup \set{i + 1, \ldots, n} \right)
\]
and
\[
B \cap \set{2, \ldots, i}
\]
are classes of $\sigma$.

\br
For $\pi \in NC_0(n)$, define the partition $\rho(\pi) \in \NC_0(\set{2, 3, \ldots, n, \bar{1}})$ (where $\bar{1}$ is identified with $n+1$) by
\begin{align*}
& i \stackrel{\pi}{\sim} j \Leftrightarrow i \stackrel{\rho(\pi)}{\sim} j \text{ for } i, j \neq 1, \\
& 1 \stackrel{\pi}{\sim} j \Leftrightarrow \bar{1} \stackrel{\rho(\pi)}{\sim} j.
\end{align*}
Clearly
\[
\rho(\NC_0'(n)) = \set{\pi \in \NC_0(\set{2, \ldots, n, \bar{1}}) | n \stackrel{\pi}{\sim} \bar{1}},
\]
and
\[
\rho(A(n)) = \set{\pi \in \NC_0(\set{2, \ldots, n, \bar{1}}) | 2 \stackrel{\pi}{\sim} n \stackrel{\pi}{\sim} \bar{1}}.
\]
In particular, $\rho(A(n)) \subset \NC_0'(\set{2, \ldots, n, \bar{1}})$. For any partition $\sigma \in \NC_0'(\set{2, \ldots, n, \bar{1}}) \backslash \rho(A(n))$, define the partition $r(\sigma) \in \rho(A(n))$ as follows: if $\bar{1} \in B \in \sigma$ and $n \in C \in \sigma$, let $r(\sigma)$ be the partition with the same classes as $\sigma$ except that $B \cup C$ is a class of $r(\sigma)$. Conversely, any such $\sigma$ can be obtained by starting with $\pi \in \rho(A(n))$ with the (unique) outer class $B$, choosing $i \in B$, $2 < i < n$ (if it exists), and taking $\sigma$ to have the same classes as $\pi$ except
\[
B \cap \left( \set{2, \ldots, i} \cup \set{\bar{1}} \right)
\]
and
\[
B \cap \set{i+1, \ldots, n}
\]
are classes of $\sigma$.

\br
We will use the usual commutator notation $[T,T'] = T T' - T' T$. Using equation~\eqref{Trace1},
\[
\begin{split}
& \ip{e_{u(1)}}{T_{u(2)} T_{u(3)} T_{u(4)} \ldots T_{u(n-2)} T_{u(n-1)} e_{u(n)}} \\
&\quad = \ip{e_{u(2)}}{T_{u(1)} T_{u(3)} T_{u(4)} \ldots T_{u(n-2)} T_{u(n-1)} e_{u(n)}} \\
&\quad = \ip{e_{u(2)}}{\bigl[T_{u(1)}, T_{u(3)}\bigr] T_{u(4)} \ldots T_{u(n-1)} e_{u(n)}} + \ldots \\
&\qquad + \ip{e_{u(2)}}{T_{u(3)} T_{u(4)} \ldots T_{u(n-2)} \bigl[T_{u(1)}, T_{u(n-1)}\bigr] e_{u(n)}} \\
&\qquad + \ip{e_{u(2)}}{T_{u(3)} T_{u(4)} \ldots T_{u(n-1)} T_{u(n)} e_{u(1)}}
\end{split}
\]
Now using equation~\eqref{Commutator}, $[T_i, T_j] = c (E_{ij} - E_{ji})$,
\[
\begin{split}
& \ip{e_{u(1)}}{T_{u(2)} T_{u(3)} T_{u(4)} \ldots T_{u(n-2)} T_{u(n-1)} e_{u(n)}} \\
&\quad = c \ip{e_{u(2)}}{e_{u(1)}} \ip{e_{u(3)}}{T_{u(4)} \ldots T_{u(n-1)} e_{u(n)}} \\
&\qquad + c \ip{e_{u(2)}}{T_{u(3)} e_{u(1)}} \ip{e_{u(4)}}{T_{u(5)} \ldots T_{u(n-1)} e_{u(n)}} \\
&\qquad + c \ip{e_{u(2)}}{T_{u(3)} \ldots T_{u(n-2)} e_{u(1)}} \ip{e_{u(n-1)}}{e_{u(n)}} + \ldots \\
&\qquad + \ip{e_{u(2)}}{T_{u(3)} \ldots T_{u(n)} e_{u(1)}} \\
&\qquad - c \ip{e_{u(2)}}{e_{u(3)}} \ip{e_{u(1)}}{T_{u(4)} \ldots T_{u(n-1)} e_{u(n)}} \\
&\qquad - c \ip{e_{u(2)}}{T_{u(3)} e_{u(4)}} \ip{e_{u(1)}}{T_{u(5)} \ldots T_{u(n-1)} e_{u(n)}} - \ldots \\
&\qquad - c \ip{e_{u(2)}}{T_{u(3)} \ldots T_{u(n-2)} e_{u(n-1)}} \ip{e_{u(1)}}{e_{u(n)}}
\end{split}
\]
The left-hand-side of the preceding equation is equal to $\Theta(\mb{\hat{1}_n}; \set{1, \ldots, n}, \vec{u})$, where $\mb{\hat{1}_n} \in \NC_0'(n)$ is the partition with a single class. Similarly, using Lemma~\ref{Lemma:Last}, the preceding equation itself states that
\[
\begin{split}
\Theta(\mb{\hat{1}_n}; \set{1, \ldots, n}, \vec{u})
& = \Theta(\rho(\mb{\hat{1}_n}); \set{2, \ldots, n, \bar{1}}, \rho(\vec{u})) \\
&\quad + \sum_{i=2}^{n-2} \Theta(\bigl(\set{2, \ldots, i, \bar{1}}, \set{i+1, \ldots, n} \bigr) \set{2, \ldots, n, \bar{1}}, \rho(\vec{u})) \\
&\quad - \sum_{i=3}^{n-1} \Theta(\bigl(\set{2, \ldots, i}, \set{1, i+1, \ldots, n}\bigr); \set{1, \ldots, n}, \vec{u}),
\end{split}
\]
where for a multi-index $\vec{u} = ((u(1), \ldots, u(n))$, we denote by $\rho(\vec{u}) = (u(2), \ldots, u(n), u(1))$ the multi-index on $\set{2, \ldots, n, \bar{1}}$. Using the descriptions of $l(\sigma)$, $r(\sigma)$ at the beginning of the proof, this in turn equals to
\begin{equation}
\label{Single-class}
\begin{split}
& = \Theta(\rho(\mb{\hat{1}_n}); \set{2, \ldots, n, \bar{1}}, \rho(\vec{u})) \\
&\quad + \sum_{\tau: r(\tau) = \rho(\mb{\hat{1}_n})} \Theta(\tau; \set{2, \ldots, n, \bar{1}}, \rho(\vec{u})) \\
&\quad - \sum_{\sigma: l(\sigma) = \mb{\hat{1}_n}} \Theta(\sigma); \set{1, \ldots, n}, \vec{u}).
\end{split}
\end{equation}
Using Lemma~\ref{Lemma:Last} again and equation~\eqref{Single-class} applied to the class of $\pi$ containing $1$, we conclude that for $\pi \in A(n)$,
\[
\begin{split}
c^{-(\abs{\pi} - 1)} \Theta(\pi; \set{1, \ldots, n}, \vec{u})
& = c^{-(\abs{\pi} - 1)} \Theta(\rho(\pi); \set{2, \ldots, n, \bar{1}}, \rho(\vec{u})) \\
&\quad + c \sum_{\tau: r(\tau) = \rho(\pi)} c^{-(\abs{\tau} - 1)} \Theta(\tau; \set{2, \ldots, n, \bar{1}}, \rho(\vec{u})) \\
&\quad - c \sum_{\sigma: l(\sigma) = \pi} c^{-(\abs{\sigma} - 1)} \Theta(\sigma; \set{1, \ldots, n}, \vec{u}),
\end{split}
\]
or, since $\abs{\rho(\pi)} = \abs{\pi}$ and $\abs{\tau} = \abs{\sigma} = \abs{\pi} + 1$,
\[
\begin{split}
\Theta(\pi; \set{1, \ldots, n}, \vec{u}) & = \Theta(\rho(\pi); \set{2, \ldots, n, \bar{1}}, \rho(\vec{u})) \\
& + \sum_{\tau: r(\tau) = \rho(\pi)} \Theta(\tau; \set{2, \ldots, n, \bar{1}}, \rho(\vec{u})) \\
& - \sum_{\sigma: l(\sigma) = \pi} \Theta(\sigma); \set{1, \ldots, n}, \vec{u}),
\end{split}
\]
Therefore
\[
\begin{split}
& \ip{e_{u(1)}}{S_{u(2)} \ldots S_{u(n-1)} e_{u(n)}} \\
&\quad = \sum_{\pi \in \NC_0'(n)} \Theta(\pi; \set{1, \ldots, n}, \vec{u}) \\
&\quad = \sum_{\pi \in A(n)} \Bigl( \Theta(\pi; \set{1, \ldots, n}, \vec{u}) + \sum_{\sigma: l(\sigma) = \pi} \Theta(\sigma; \set{1, \ldots, n}, \vec{u}) \Bigr) \\
&\quad = \sum_{\pi \in A(n)} \Bigl( \Theta(\rho(\pi); \set{2, \ldots, n, \bar{1}}, \rho(\vec{u})) + \sum_{\tau: r(\tau) = \rho(\pi)} \Theta(\tau; \set{2, \ldots, n, \bar{1}}, \rho(\vec{u})) \Bigr) \\
&\quad = \sum_{\tau \in \rho(\NC_0'(\set{2, \ldots, n, \bar{1}}))} \Theta(\tau; \set{2, \ldots, n, \bar{1}}, \rho(\vec{u})) \\
&\quad = \ip{e_{u(2)}}{S_{u(3)} \ldots S_{u(n)} e_{u(1)}}.
\end{split}
\]
Using Theorem~\ref{Thm:Cumulants}, we conclude that the free cumulant functional $R_\phi$, and so $\phi$ itself, is tracial.
\end{proof}

\begin{Ex}
For $d=2$, one can take
\[
T_1 =
\begin{pmatrix}
c+1 & 0 \\
0 & 1
\end{pmatrix},
\qquad
T_2 =
\begin{pmatrix}
0 & 1 \\
1 & 0
\end{pmatrix}.
\]
For $d=3$, one can take
\[
T_1 =
\begin{pmatrix}
0 & c & 0 \\
c & 0 & 0 \\
0 & 0 & 0
\end{pmatrix},
\qquad
T_2 =
\begin{pmatrix}
c & 0 & 0 \\
0 & c+1 & 0 \\
0 & 0 & 1
\end{pmatrix},
\qquad
T_3 =
\begin{pmatrix}
0 & 0 & 0 \\
0 & 0 & 1 \\
0 & 1 & 0
\end{pmatrix}.
\]
For $d=4$, one can take
\[
T_1 =
\begin{pmatrix}
c & 0 & 1 & 0 \\
0 & 0 & 0 & 0 \\
1 & 0 & 0 & 0 \\
0 & 0 & 0 & 1
\end{pmatrix},
\;
T_2 =
\begin{pmatrix}
0 & 0 & 0 & 0 \\
0 & 0 & c & 0 \\
0 & c & 0 & 0 \\
0 & 0 & 0 & 0
\end{pmatrix},
\;
T_3 =
\begin{pmatrix}
1 & 0 & 0 & 0 \\
0 & c & 0 & 0 \\
0 & 0 & c+1 & 0 \\
0 & 0 & 0 & 1
\end{pmatrix},
\;
T_4 =
\begin{pmatrix}
0 & 0 & 0 & 1 \\
0 & 0 & 0 & 0 \\
0 & 0 & 0 & 1 \\
1 & 0 & 1 & 0
\end{pmatrix}.
\]
\end{Ex}

\begin{Remark}
If $C_{ij} = c$ for all $i, j$, the corresponding Fock space is an interacting Fock space in the sense of \cite{AccBozGaussianization}. If $C_{ij} = c$ and in addition all $T_i = 0$, the von Neumann algebras $W^\ast(X_1, \ldots, X_d)$ were described in \cite{Ricard-t-Gaussian}.
\end{Remark}


\providecommand{\bysame}{\leavevmode\hbox to3em{\hrulefill}\thinspace}
\providecommand{\MR}{\relax\ifhmode\unskip\space\fi MR }
\providecommand{\MRhref}[2]{%
  \href{http://www.ams.org/mathscinet-getitem?mr=#1}{#2}
}
\providecommand{\href}[2]{#2}

\end{document}